\documentclass{article}
\usepackage{amsmath, amsfonts, amscd}

\newtheorem{theorem}{Theorem}[section]
\newtheorem{lemma}{Lemma}[section] 
\newtheorem{definition}{Definition}[section] 
\newtheorem{corollary}[theorem]{Corollary}
\newcommand{\problem}[1]{\\[2 mm]\textbf{Problem: }{\em #1 }} 
\numberwithin{equation}{section}
\newcommand{\h}[1]{\hspace{#1 truein}}  
\newcommand{\sss}{\scriptscriptstyle} 

\newcommand{\freccia}{\longrightarrow} 
\newcommand{\xfreccia}[1]{\xrightarrow{\ \ #1\ \ }} 
\newcommand{\xfrecciad}[1]{\xrightarrow{\displaystyle\ \ #1\ \ }}

\newcommand{\eps}{\varepsilon} 
\newcommand{\field}[1]{\mathbb{#1}}
\newcommand{\R}{\field{R}}                        
\newcommand{\ER}{{^\bullet\R}}                    
\newcommand{\N}{\field{N}}                        
\newcommand{\diff}[1]{\,\hbox{\rm d}#1}               

\newcommand{\st}[1]{{^\circ #1}} 
\newcommand{\Top}[1]{{\mbox{\Large$\tau$}}_{\sss{#1}}} 
\newcommand{\In}[1]{\in_{_{\sss{#1}}}} 
\newcommand{\InUp}[1]{\in^{\sss{#1}}} 
\newcommand{\Nil}{\mathcal{N}} 
\newcommand{\No}[1]{\mathcal{N}_{#1}} 
\newcommand{\C}{\mathcal{C}} 
\newcommand{\D}{\mathcal{D}} 
\newcommand{\F}{\mathcal{F}} 
\newcommand{\Ds}[1]{\mathcal{D}_{\sss #1}} 
\newcommand{\Fs}[1]{\mathcal{F}_{\sss #1}} 
\newcommand{\ext}[1]{{}^\bullet #1} 
\DeclareMathOperator{\Lip}{Lip} 

\DeclareMathOperator*{\Man}{{\bf Man}\mbox{${}^n$}}
\DeclareMathOperator*{\ManInfty}{{\bf Man}}
\DeclareMathOperator{\Set}{{\bf Set}}
\DeclareMathOperator*{\ORn}{{\bf O}\mbox{$\R^n$}}
\DeclareMathOperator*{\ORInfty}{{\bf O}\mbox{$\R^\infty$}}
\newcommand{\Cn}{\boldsymbol{\mathcal{C}}^n} 
\newcommand{\CInfty}{\boldsymbol{\mathcal{C}}^{\infty}} 
\newcommand{\ECn}{\ext{\Cn}} 
\newcommand{\ECInfty}{\ext{\CInfty}} 
\DeclareMathOperator{\SERn}{{\bf S}\mbox{$\ER^n$}} 

\newcommand{\Lim}[2]{\lim_{t \to 0}\frac{{#1}(t)-{#2}(t)}{t}} 
\newcommand{\Limup}[1]{\lim_{t \to 0}\frac{#1}{t}} 
\newcommand{\Limsup}[2]{\limsup_{t \to 0}\left|\frac{{#1}(t)-{#2}(t)}{t}\right|}
\newcommand{\Limsupup}[1]{\limsup_{t \to 0}\left|\frac{#1}{t}\right|}

\newcommand{\qed}{\\\phantom{|}\hfill \vrule height 1.5ex width 1.3ex depth-.2ex \\}
\newcommand{\qedNoNewLine}{\\\phantom{|}\hfill \vrule height 1.5ex width 1.3ex depth-.2ex}
\newcommand{\e}[1]{\text{\h{.1}\  #1 \h{.1}}} 
\newcommand{\ee}[1]{\text{\h{.3}\  #1 \h{.3}}} 
\newcommand{\ptind}{\displaystyle \mathop {\ldots\ldots\,}} 
\newcommand{\pti}{:\;\;\;} 
\newcommand{\then}{\quad \Longrightarrow \quad}
\newcommand{\DIff}{ \quad\;\; :\!\iff \quad } 
\newcommand{\Proof}{\ \ \emph{Proof: }}

\title{\textbf{Infinitesimal Differential Geometry\\ }}
\author{Paolo Giordano\\
\\
\emph{Inst. f. Angew. Mathematik, Univ. Bonn, Germany}\\
\emph{Accademia di Mendrisio, Univ. Svizzera Italiana, Switzerland}\\
\emph{e-mail: pgiordano@arch.unisi.ch}}
\date{}
\begin{document}
\maketitle
\thispagestyle{empty}
\begin{center}Shortened title: Infinitesimal Diff. Geom.\end{center}
\renewcommand{\thefootnote}{}
\footnotetext[1]{Mathematical Subject Classification: 58D15, 58B10, 58A05, 26E15.\\
Key words: Spaces of mappings; Nilpotent infinitesimals; Differential manifolds, foundations. This research was
supported through a DAAD (German Academic Exchange Service) and a European Community Marie Curie fellowships.\\}
\begin{abstract}Using standard analysis only, we present an extension $\ER$ of the real field containing nilpotent
infinitesimals. On the one hand we want to present a very simple setting to formalize infinitesimal methods in
Diffe\-ren\-tial Geometry, Analysis and Physics. On the other hand we want to show that these infinitesimals may be
also useful in infinite dimensional Differential Geometry, e.g. to study spaces of mappings. We define a full
embedding of the category $\Man$ of finite dimensional $\Cn$ manifolds in a cartesian closed category. In it we have a
functor $\ext(-)$ which extends these spaces adding new infinitesimal points and with values in another full cartesian
closed embedding of $\Man$. We present a first development of Differential Geometry using these infinitesimals.
\end{abstract}
\section{The ring of standard infinitesimals}
\subsection{Introduction}
Frequently in Physics it is possible to find informal calculations like
\[
    \frac{1}{\sqrt{1 - {\displaystyle \frac{v^2}{c^2}}}} = 1 + \frac{v^2}{2 c^2} \qquad \qquad
    \sqrt{1-h_{44}(x)} = 1 - \frac{1}{2} h_{44}(x)
\]
with explicit use of infinitesimals $v / c \ll 1$ or $h_{44}(x) \ll 1$ such that e.g. $h_{44}(x)^2 =0$. In fact using
this type of infinitesimals we can write an \emph{equality}, in some infinitesimal neighborhood, between a
smooth function and its tangent straight line, or, in other words, a Taylor formula without remainder.\\
Informal methods based on actual infinitesimals are sometimes used in Differential Geometry too. Some classical
examples are the following: a tangent vector is an infinitesimal arc of curve traced on the manifold and the sum of
tangent vectors is made using infinitesimal parallelograms; tangent vectors to the tangent bundle are infinitesimal
squares on the manifold; a vector field is sometimes intuitively treated as an ``infinitesimal transformation'' of the
space into itself and the Lie brackets of two vector fields as the commutator of the corresponding infinitesimal
transformations.\\
There are obviously many possibilities to formalize this kind of intuitive reasonings, obtaining a more or
less good dialectic between informal and formal thinking (see e.g. \cite{Lav, Kr-Mi} and references therein).\\
We want to show how it is possible to extend the real field adding nilpotent infinitesimals by means of a very simple
construction completely inside ``standard mathematics'' (with this we mean that the formal control necessary to work
in our setting is very less strong than that necessary both in Non-Standard Analysis \cite{Al-Fe-Ho-Li} and Synthetic
Differential Geometry \cite{Lav}). To define the extension $\ER \supset \R$ we shall use elementary analysis only.\\
The usefulness of this extension can be glimpsed saying e.g. that using $\ER$ it is possible to write in a completely
rigorous way that a smooth function is equal to its tangent straight line in a first order neighborhood, to use
infinitesimal Taylor formulas without remainder, to define a tangent vector as an infinitesimal curve and sum them
using infinitesimal parallelograms, to see a vector field as an infinitesimal transformation, hence, to come to the
point, to formalize many non-rigorous methods used in Physics and Geometry. This is important both for didactical
reasons and because it was by means of these methods that mathematicians like S. Lie and E. Cartan were originally
conducted to construct important concepts of Differential Geometry.\\
We can use the infinitesimals of $\ER$ not only as a good language to reformulate well-known results, but also as a
very useful tool to construct, in a simple and meaningful way, a Differential Geometry in classical
infinite-dimensional objects like {\bf Man}$(M,N)$ the space of all the $\C^{\infty}$ mapping between two mani\-folds
$M$, $N$. Here with ``simple and meaningful'' we mean the idea to work directly on the geometric object in an
intrinsic way without being forced to use charts, but using infinitesimal points (see \cite{Lav}). Some important
examples of spaces of mappings used in applications are the space of configurations of a continuum body, groups of
diffeomorphisms used in hydrodynamics, magnetohydrodynamics, electromagnetism, plasma dynamics and paths spaces for
calculus of variations (see \cite{Kr-Mi, Ma-Ra} and references therein). Interesting applications in classical field
theories can also be found in \cite{Ab-Ma}.\\
A complete and powerful setting for this kind of problems, but without the use of infinitesimals, can be found in
\cite{Fr-Kr, Kr-Mi}. The construction of our categories takes a strong inspiration from this works and from
\cite{Che}. The author hope that this work could also serve to introduce infinitesimal methods in the convenient
setting of \cite{Kr-Mi}. The most complete use of infinitesimals in Differential Geometry can be found in \cite{Koc,
Lav, Mo-Re}, whose setting is incompatible with classical logic and admits models in intuitionistic logic only. The
infinitesimals methods formalized in this work are strongly influenced by \cite{Lav, Koc}.

We start from the idea that a smooth ($\C^\infty$) function $f: \R \freccia \R$ is actually \emph{equal} to its
tangent straight line in the first order neighborhood e.g. of the point $x=0$, that is
\begin{equation}
\forall h \in D \pti  f(h) = f(0) + h \cdot f^\prime(0)                   \label{eq: FunctionEqualTangent}
\end{equation}
where $D$ is the subset of $\ER$ which defines the above-mentioned neighborhood of $x=0$. The previous (\ref{eq:
FunctionEqualTangent}) can be seen as a first-order Taylor formula without remainder because intuitively we think that
$h^2=0$ for any $h \in D$. These almost trivial considerations lead us to understand many things: $\ER$ must
necessarily be a ring and not a field; moreover we will surely have some limitation in the extension of some function
from $\R$ to $\ER$, e.g. the square root. But we are also led to ask if (\ref{eq: FunctionEqualTangent}) uniquely
determines the derivative $f^\prime(0)$: because even if it is true that we cannot simplify by $h$, we know that the
polynomial coefficient of a Taylor formula are unique in classical analysis. In fact we will prove that
\begin{equation}
\exists !\, m \in \R :\ \forall h \in D \pti  f(h) = f(0) + h \cdot m,                       \label{eq: IdeaDF}
\end{equation}
that is the slope of the tangent is uniquely determined in case it is an ordinary real number.\\
If we try to construct a model for (\ref{eq: IdeaDF}) a natural idea is to think our new numbers as equivalence
classes $[h]$ of usual functions $h : \R \freccia \R$. In such a way we can hope both to include the real field using
classes generated by constant functions, and that the class generated by $h(t) = t$ could be a first order
infinitesimal number. To understand how to define this equivalence relation we can see (\ref{eq:
FunctionEqualTangent}) in the following sense:
\begin{equation}
  f(h(t)) \sim f(0) + h(t) \cdot f^\prime(0).         \label{eq: IdeaFunctions}
\end{equation}
If we think $h(t)$ ``sufficiently similar to $t$'', we can define $\sim$ so that (\ref{eq: IdeaFunctions}) is
equivalent to
\[
  \Limup{f(h(t)) - f(0) - h(t) \cdot f^\prime(0)}=0,
\]
that is
\begin{equation}
  x\sim y \DIff \Lim{x}{y}=0.               \label{eq: IdeaRelEq}
\end{equation}
In this way (\ref{eq: IdeaFunctions}) is very near to the definition of differentiability for $f$ at 0.\\
It is important to note that, because of l'H\^opital's theorems
\[
  \C^1(\R,\R)/\!\sim
  \,\,\,\simeq\,\,
  \R[x]/\langle x^2 \rangle
\]
that is the usual tangent bundle of $\R$ and thus we obtain nothing new. It is not easy to understand what set of
functions we have to choose for $x$, $y$ in (\ref{eq: IdeaRelEq}) so as to obtain a non trivial structure. The first
idea is to take continuous functions at $t=0$ so that e.g. $h_k(t) = |t|^{1/k}$ is a $k$th order nilpotent
infinitesimal; for almost all the results presented in this article continuous functions at $t=0$ work well, but only
in proving the non-trivial property
\begin{equation}
    \label{thm: UniquenessIncrementalRatio}
    \left(\forall x \in \ER \pti x \cdot f(x) =0 \right) \then \forall x \in \ER \pti f(x)=0
\end{equation}
(here $f : \ER \freccia \ER$ is a smooth function, in a sense we shall precise after) we will see that it doesn't
suffice to take continuous functions at $t=0$. The previous property (\ref{thm: UniquenessIncrementalRatio}) is
useful to prove the uniqueness of smooth incremental ratios, hence to define the derivative $f^\prime: \ER \freccia
\ER$ for a smooth function $f : \ER \freccia \ER$ which, generally speaking, is not the extension to $\ER$ of an
ordinary function defined on $\R$ (e.g. the function used for the small oscillations of the pendulum $t \mapsto \sin(h
\cdot t)$, where $h \in \ER \setminus \R$). To prove (\ref{thm: UniquenessIncrementalRatio}) the following functions
turned out to be very useful:
\begin{definition}
If $x: \R \freccia \R$, then we say that $x$ \emph{is nilpotent} iff $|x(t) -x(0)|^k={\rm o}(t)$ for some $k \in \N$.
$\Nil$ will be the set of all the nilpotent functions.
\end{definition}
E.g. any Holder function $|x(t) - x(s)| \le c \cdot |t - s|^\alpha$ ($\alpha > 0$) is nilpotent. Hence we now define
\begin{definition}
Let $x,y \in \Nil$, then we say $x \sim y$ iff
\[
    x(t)=y(t)+\text{\rm o}(t)\e{for}t \to 0.                    \label{def: sim}
\]
The quotient $\Nil / \sim$ will be indicated with $\ER$ and called ``\emph{the ring of standard infinitesimals}''. Its
elements $x \in \ER$ will be called ``\emph{extended reals}''. We can read $\ER$ either as ``dot R'' or ``extended
R''.
\end{definition}
E.g. the previous $h_k(t)=|t|^{1/k}$ is not equivalent to zero but its $k+1$-th power is equivalent to zero, thus it
is a nilpotent infinitesimal. Because it is also an ordinary infinitesimal function for $t \to 0$ this motivates the
name ``ring of \emph{standard} infinitesimals''. $\Nil$ is close with respect to pointwise sum and product of
functions. For the product it suffices to write $x\cdot y -x(0) \cdot y(0) = x \cdot [y - y(0)] + y(0) \cdot
[x-x(0)]$. The case of the sum follows from the subsequent equalities (where we use $x_t:=x(t)$, $u := x-x_0$ and
$v:=y-y_0$):
\begin{gather*}
u^k \sim 0 \sim v^k\\
(u+v)^k=\sum_{i=0}^k \binom{k}{i} u^i \cdot v^{k-i}\\
\forall i=0, \dots, k \pti \frac{u_t^i \cdot v_t^{k-i}}{t} = \frac{\left(u_t^k\right)^{\frac{i}{k}} \cdot
\left(v_t^k\right)^{\frac{k-i}{k}}}{t^{\frac{i}{k}} \cdot t^{\frac{k-i}{k}}} =
\left(\frac{u_t^k}{t}\right)^{\frac{i}{k}} \cdot \left(\frac{v_t^k}{t}\right)^{\frac{k-i}{k}}.
\end{gather*}
Obviously $\sim$ is a congruence relation with respect to pointwise operations hence $\ER$ is a commutative ring.\\
Where it will be useful to simplify notations we will write ``$x = y$ in $\ER$'' instead of $x \sim y$, and we
will talk directly about the elements of $\Nil$ instead of their equivalence classes; for example we can say that
$x=y$ in $\ER$ and $z=w$ in $\ER$ imply $x+z=y+w$ in $\ER$.\\
The immersion of $\R$ in $\ER$ is $r \longmapsto \hat r$ defined by $\hat r(t) := r$, and in the sequel we will always
identify $\hat{\R}$ with $\R$. Conversely if $x \in \ER$ then is well defined and meaningful the {\em standard part
map\/} $\st(-) : x \in \ER \longmapsto \st x = x(0) \in \R$ which evaluates each extended real in $0$.
\subsection{The ideal of first order infinitesimals}
If we want that $f(h(t)) \sim f(0) + h(t) \cdot f^\prime(0)$ then from Taylor formula we obtain
\begin{equation}
  \Limup{f(h(t)) - f(0) - h(t) \cdot f^\prime(0)} = \Limup{h(t)}\cdot \sigma(t)                 \label{eq: TaylorIdea}
\end{equation}
with $\sigma(t) \to 0$ for $t \to 0$. This suggests us to define $D$ using the condition
\[
  \Limsupup{h(t)} < +\infty.
\]
Generally we will write $x \approx y$ for
\[
    \Limsup{x}{y} < +\infty
\]
and we will say that \emph{$x$ is close to $y$}. We obtain a well-defined congruence on $\ER$ that coincides with
equality on $\R$.
\begin{definition}
\[
  D := \{h \in \ER\,|\, h \approx 0 \}.
\]
The elements of $D$ are called {\em first order infinitesimals}.
\end{definition}
Thus we have $x \approx y$ iff $x(t)=y(t)+{\rm O}(t)$ for $t \to 0$. For example if $r,s \in \R$, then $h(t) := r
\cdot |t|$ if $t \ge0$ and $h(t) := s \cdot |t|$ if $t \le 0$ is a first order infinitesimal; another one is $h(t) :=
r \cdot t \cdot \sin(1/t)$, and obviously $h(t) := t$ and in general any $\C^1$ infinitesimal function at $t=0$.
Conversely, if $\alpha \in (1/2,1)$, then $x(t) := |t|^\alpha$ is not an element of $D$ but note that $x^2=0$ in
$\ER$.
\begin{theorem}
\label{thm: DIdeal}
$D$ is an ideal of $\ER$, and
\[
  \forall h \in D \pti h^2=0.
\]
\end{theorem}
\Proof It follows from elementary properties of $\limsup$; for example the inequalities
$$
  \Limsup{h}{k} \le \Limsupup{h(t)}+\Limsupup{k(t)}<+\infty
$$
$$
  \Limsupup{x(t)\cdot h(t)}\le |x(0)| \cdot \Limsupup{h(t)}<+\infty
$$
prove that $D$ is an ideal, and the following
$$
           0 \le \liminf_{t \to 0}\left|\frac{h(t)^2}{t}\right| \le
             \Limsupup{h(t)^2} \le
             |h(0)|\cdot\Limsupup{h(t)} =0
$$
prove that every element of $D$ has square equal to zero. \qed
Another interesting ideal is $D_k := \{ h \in \ER \,|\,h^k \in D \}$ for $k \in \N_{\sss >0}$: this follows from
Newton's formula and the equality
\begin{equation}
\label{eq: NewtonPowers}
    \frac{h(t)^i \cdot u(t)^{k-i}}
         {t}
          =
    \left[
            \frac{h(t)^{k}}
                 {t}
   \right]^\frac{i}{k} \cdot  \left[
                                      \frac{u(t)^{k}}
                                           {t}
                                \right]^\frac{k-i}{k}.
\end{equation}
It is also useful to define $D_0 :=\{ 0 \}$. Using an idea similar to (\ref{eq: NewtonPowers}) and taking $h_k \in
D_{j_k}$, and $0 \le i_k$, we also have
\begin{equation}
\label{eq: PowerRules}
    \begin{split}
      {}& h_1^{i_1} \cdot \ldots \cdot h_n^{i_n} =0        \ee{if} \sum_{k=1}^n\frac{i_k}{j_k} > 1 \\
      {}& h_1^{i_1} \cdot \ldots \cdot h_n^{i_n} \in D_p   \ee{if} \frac{1}{p} \le \sum_{k=1}^n\frac{i_k}{j_k} \le 1.
    \end{split}
\end{equation}
E.g. if $h \in D_3$ and $u \in D_5$ we have $h^2 u^3=0$ and $h^2 u \in D_2$. It may also useful to note that $h k =0$
if $h^2=k^2=0$ and $h i =0$ if $h \in D$ and $\st{i}=0$, that is $i$ is a generic infinitesimal. Another useful
property is expressed by the following cancellation law, which is a good substitute for the fact that $\ER$ is not a
field.
\begin{theorem}
Let $x \in \ER$ and $x \ne 0$, then
\[
  x \cdot r = x \cdot s \e{and} r,s \in \R \then r=s.                  \label{thm: CancellationLaw}
\]
\end{theorem}
\Proof We can write the hypothesis $x \cdot r = x \cdot s$ as
\[
  \Limup{x(t)}\cdot (r-s) = 0 = |r-s| \cdot \Limsupup{x(t)},
\]
but the $\Limsupup{x(t)} \ne 0$ because $x \ne 0$, and hence $r = s$. \qed
Obviously this law is not true if $r$, $s$ are generic extended reals. Finally it is also easy to prove that
$x \in \ER$ is invertible iff $\st{x}\ne 0$.
\subsection{Extension of functions}
Before considering the proof of (\ref{eq: IdeaDF}) we have to understand how to extend a given function $f : \R
\freccia \R$ to a certain $\ext f : \ER \freccia \ER$. First of all we can define $\ext A$ for $A \subseteq \R^k$
exactly as we defined $\ER$: it is sufficient to consider the set $\No{A}$ of all the nilpotent functions $ x : \R
\freccia A$ (that is such that $||x_t - x_0 ||^k = \text{o}(t)$ for some $k \in \N$, where $||-||$ is the norm in
$\R^k$) with values in $A$; afterward we take the quotient with respect to the analogous of the relation $\sim$
defined in Def. \ref{def: sim}. We shall give further the general definition of the extension functor $\ext(-)$, here
we only want to examine some elementary properties of the ring $\ER$.
\begin{definition}
Let $A$ be a subset of $\R^k$, $f : A \freccia \R$ and $x \in \ext A$ then we define
\[
 \ext f (x) := f \circ x.
\]
\end{definition}
This definition is well posed if $f$ is locally lipschitzian; in fact if $x=y$ in $\ER$ then $x_0 = y_0$ and so for
some $\delta, K >0$ we have
\begin{equation}
    \forall t \in (-\delta, \delta) \pti \|f(x_t) - f(y_t) \| \le K \cdot \| x_t - y_t \|;          \label{eq: LocallyLip}
\end{equation}
hence for $t \in (-\delta, \delta)$ we have
\begin{align*}
           0&\le\displaystyle \Limsupup{f(x_t) - f(y_t)} \le   \\*[3 mm]
          {}&\le\displaystyle K \cdot \Limsup{x}{y} = 0.
\end{align*}
Note also that (\ref{eq: LocallyLip}) implies $f(x) \in \Nil$ if $x \in \No{A}$. In the sequel $\Lip(A,B)$ will be the
set of all the locally lipschitzian functions defined in $A$ and with values in $B$. The function $\ext f$ is an
extension of $f$, that is
\[
  \ext f(r) = f(r) \quad \rm{in} \quad \ER \quad \rm{for} \quad r \in \R,
\]
thus we can still use the symbol $f(x)$ both for $x \in \ER$ and $x \in \R$ without confusion.\\
In the following theorem $I_0 := \{ h \in \ER \,|\, \st{h}=0 \}$ will be the set of all the infinitesimals of $\ER$.
\begin{theorem}
    Let $A$ be an open set in $\R$ and $x \in A$, then $x + h \in \ext{A}$ for every $h \in I_0$.     \label{thm: DandOpen}
\end{theorem}
It is necessary to give some explanation to understand the statement of this theorem. In fact $\ext{A}=\No{A} / \sim$,
thus we don't have $\ext{A} \subseteq \ER$ if $A \subseteq \R$ (any equivalence relation $[x]_A \in \ext{A}$ is made
of functions $x : \R \freccia A$ only, whereas $[x]_{\R} \in \ER$ is made of functions $x : \R \freccia \R$). In spite
of all that there is obviously a natural injection $i : \ext{A} \freccia \ER$. In fact $[x]_A = \left\{ y \in \No{A}
\,|\, x \sim y \right\}$ and so $x \in \Nil=\No{\R}$ and we can define $i([x]_A):=[x]_{\R}$. This map is well defined
and injective, essentially because the definition of $\sim$ doesn't depend on $A$. Using $i : \ext{A} \freccia \ER$ we
can identify $\ext{A}$ with a subset of $\ER$ if it is clear from the context the superset we are considering (in this
case $\R \supseteq A$); the statement of the previous theorem use this identification.\\
\Proof We have to prove that $[x+h]_{\R} \in i(\ext{A})$. Because $h \in I_0$ we have that $x + h_t \in A$ for $t$
sufficiently small $t \in (-\delta, \delta)$ and thus there exists $y : \R \freccia A$ such that $y_t = x + h_t$ for
$t \in (-\delta, \delta)$. Hence, directly from the definition of $\sim$, $i([y]_A)=[y]_{\R}=[x+h]_{\R}$.\qed In
conclusion of this section we enunciate the following useful elementary transfer theorem for equalities, whose proof
follows directly from the previous definitions:
\begin{theorem}
Let $A \subseteq \R^k$, and $\tau, \sigma \in \Lip(A,\R)$. Then it results
\[
  \forall x \in {\ext A} \pti  \ext\tau(x) = \ext\sigma(x)                              \label{thm: transfer=}
\]
iff
\[
  \forall r \in A \pti \tau(r) = \sigma(r).
\]
\end{theorem}
\subsection{The derivation formula}
Now we will prove the formula (\ref{eq: IdeaDF}), which we will call {\em derivation formula}. It is natural to expect
that it will be equivalent to the usual differentiability of a function, in fact we have
\begin{theorem}
    \label{thm: DerivationFormula}
    Let $A$ be an open set in $\R$, $x \in A$ and $f \in \Lip(A;\R)$, then the following
    conditions are equivalent:
    \begin{enumerate}
        \item \ \ $f$ is differentiable at x
        \item \ \ $\exists!\,m \in \R :\; \forall h \in D \pti f(x+h) = f(x) + h \cdot m$.
    \end{enumerate}
    In this case we have $m = f^\prime(x)$, where $f^\prime(x)$ is the usual derivative of $f$ at $x$.
\end{theorem}
Note that $m = f^\prime(x) \in \R$, i.e. the slope is an usual real number and that we can use the previous formula
with standard real numbers $x$ only, and not with a generic $x \in \ER$, but we shall remove this limitation in a
subsequent section. In other words we can say that this formula allows us to differentiate the
usual differentiable functions using a language with infinitesimal numbers and to obtain from this an ordinary function.\\
\Proof 1) $\Rightarrow$ 2): First of all note that because of Theorem \ref{thm: DandOpen} we can consider $f(x+h)$ for
any $h \in D$. Now let $m := f^\prime(x)$ and $h \in D$, i.e. $\Limsupup{h(t)} < +\infty$. For hypothesis $f$ is
differentiable in $x$, hence we can find a function $\sigma : (A-x)\freccia \R$ such that
\[
    \forall u \in A-x \pti f(x+u) = f(x) + u \cdot m + u \cdot \sigma(u)
\]
\[
    \lim_{u \to 0} \sigma(u)=\sigma(0)=0.
\]
Therefore
\begin{align*}
    \Limsupup{f(x+h_t)-f(x)-h_t \cdot m} &= \Limsupup{h_t \cdot \sigma(h_t)} \\
                                       {}&\le \sigma(h_0)\cdot\Limsupup{h_t} = 0.
\end{align*}
This proves the existence; for the uniqueness we simply use the cancellation law (Theorem \ref{thm: CancellationLaw}).\\
2) $\Rightarrow$ 1): For this implication it suffices to apply the hypothesis 2) with $h(t) := t$. \qed If we apply
this theorem to the $\C^1$ function $p(r):=\int_x^{x+r}f(t)\diff{t} $, then we obtain the following
\begin{corollary}
Let $A$ open in $\R$, $x \in A$ and $f \in \C^0(A)$. Then
\[
  \forall h \in D \pti \int_x^{x+h} f(t)\diff{t}=h \cdot f(x).
\]
Moreover $f(x) \in \R$ is uniquely determined by this equality.
\end{corollary}
 For multiple integrals we have analogous formulas;
e.g. if $h,k \in D_2$ and $h \cdot k \in D$ then
\[
    \int_{\sss [0,h]\times[0,k]}f(x,y)\diff{x}\diff{y} = h k \cdot f(0,0).
\]
With the ideal $D_k$ of the $k$th order infinitesimal numbers and a function $f \in \C^k(A)$ it is possible to prove
infinitesimal Taylor formula without any remainder
\[
  \forall h \in D_k \pti f(x+h) = \sum_{i=0}^k \frac{h^i}{i!} \cdot f^{(i)}(x)
\]
with the standard reals $f^{(i)}(x)$ uniquely determined by this formula. Another useful form of the derivation
formula is the following
\begin{theorem}
    Let $A$ open in $\R$ and $f : A \freccia \R$ be a $\C^1$ function. Let $h,k \in \ER$ be such that $h \cdot k \in D$,
    then for every $x \in A$
    \[
        k \cdot f(x + h) = k \cdot f(x) + k h \cdot f'(x)                       \label{thm: SecondDerivationFormula}
    \]
\end{theorem}

We close this section introducing a very simple notation useful to emphasize some equalities: if $h,k \in \ER$ then we
say that $\exists h/k$ iff $\exists! r \in \R \;:\; h = r \cdot k$, and obviously we indicate this $r \in \R$ with
$h/k$. Therefore we can say, e.g., that
\begin{align*}
    f^{\prime}(x) &= \frac{f(x + h) - f(x)}{h} &{} \\
                {}&                          {}& \forall h \in D_{\neq 0} \\
             f(x) &= \frac{1}{h} \cdot \int_x^{x+h} f(t)\diff{t}. &{}
\end{align*}
Moreover we can prove some natural properties of this ``ratio'', like the following one
\[
    \exists \frac{u}{v}, \frac{x}{y} \e{and} vy \neq 0 \Longrightarrow
    \frac{u}{v} + \frac{x}{y} = \frac{uy + vx}{vy}.
\]
\subsection{Order relations}
From the previous sections one can draw the conclusion that $\ER$ is essentially ``the little-oh'' calculus. If on the
one hand this is certainly true, on the other hand the extended reals give us more flexibility than this calculus:
working with $\ER$ we don't have to bother ourselves with remainders made of ``little-oh'', but we can neglect them
and use the great powerfulness of the algebraic calculus with nilpotent infinitesimals (see \cite{Lav} for many
examples which can be repeated almost equal in our setting using previous theorems). But thinking the elements of
$\ER$ as new numbers, and not simply as ``little-oh functions'', permits to treat them in a different and new way, for
example to define on them two meaningful partial order relations, the first one of which is the following.
\begin{definition}
For $x,y \in \ER$, we say that $x \succeq y$ iff we can find $z \in \ER$ such that $z=0$ in $\ER$ and
\[
  \exists \,\delta>0 \pti  \forall\, t \in (-\delta,\delta) \pti x(t) \ge y(t) + z(t).              \label{def: OrderRel}
\]
In other words let us write $\forall^{\sss{0}}\, t : \mathcal{P}(t)$ to indicate that the property $\mathcal{P}(t)$ is
true for all $t$ in some neighborhood of $t = 0$, then we can reformulate the previous definition using the
``little-oh'' language
\[
    x \succeq y \DIff \forall^{\sss{0}}\, t \pti x(t) \ge y(t) + {\rm o}(t),
\]
but note that the function ${\rm o}(t)$ depends on $x,y$. We can read $x \succeq y$ saying ``$x$ is weakly greater or
equal to $y$''.
\end{definition}
We can equivalently say that $x \succeq y$ iff we can find $x=x^\prime$ and $y=y^\prime$ in $\ER$ such that
$\forall^{\sss{0}}t:\ x^\prime_t \ge y^\prime_t$. The definition of $\succeq$ is well posed, and for example we have
that the first order infinitesimal $h(t)=\vert t \vert$ is positive but not negative. It is easy to prove that this
relation is reflexive and transitive, hence it remains to show that it is also anti-symmetric. If $x \succeq y$ and $y
\succeq x$ then we have
\begin{align*}
    {}& x(t)-y(t) \ge z_1(t) \qquad \forall t \in (-\delta_1,\delta_1) \\*[2.1 mm]
    {}& y(t)-x(t) \ge z_2(t) \qquad \forall t \in (-\delta_2,\delta_2) \\*[2.1 mm]
    {}&\displaystyle \Limup{z_1(t)} = 0 = \Limup{-z_2(t)}.
\end{align*}
Taking $\delta := \min\{\delta_1,\delta_2\}$ we obtain
\[
  \forall t \in (0,\delta) \pti
  \frac{z_1(t)}{t} \le \frac{x(t)-y(t)}{t} \le -\frac{z_2(t)}{t}
\]
\[
  \forall t \in (-\delta,0) \pti
  \frac{-z_2(t)}{t} \le \frac{x(t)-y(t)}{t} \le \frac{z_1(t)}{t}.
\]
Hence for $t \to 0$, these inequalities prove that $x=y$ in $\ER$.\\
With this relation $\ER$ becomes an ordered ring. We also observe that $\succeq$ extends the order relation in $\R$
and that it is possible to prove the cancellation law for inequality, that is if $h \in \ER$ is different from zero
and $r,s \in \R$, then from $|h| \cdot r \preceq |h| \cdot s$ we can deduce that $r \le s$.\\
We can enunciate an elementary transfer theorem for inequalities, simply substituting $=$ with $\preceq$ in Theorem
\ref{thm: transfer=}. Finally note that the usual definition of infinitesimal number as an extended real $x$ for which
$-r \prec x \prec r$ for all standard positive real number $r$ is equivalent to say that the standard part of $x$ is
zero.\\
\indent It is possible to define another meaningful partial order relation on $\ER$ saying that
\[
    x \le y \DIff x=y \quad\text{or}\quad (x \preceq y \quad\text{and}\quad y-x \text{ \ is invertible}).
\]
Some properties are better stated using $\preceq$ (e.g. elementary transfer theorem, properties of absolute value and
those about infinitesimals), whereas $\le$ is better for powers and logarithms, for topological properties and for
intervals. Actually, as we will see, a useful topology on $\ER$ is generated by the sets $\ext{U}$ for $U$ open in
$\R$; it is easy to see that if $h_t := |t \cdot \sin \frac{1}{t}|$, then $0$ is not an interior point neither in $\{x
\in \ER | -h \preceq x \preceq h \}$ nor in $\{x \in \ER | -h \preceq x \preceq 1 \}$. Therefore the above mentioned
topology is not generated by $\preceq$, whereas it is easy to check that it is generated by $\le$.\\
Once again the ring structure of $\ER$ is compatible with $\le$; the order relation between standard reals is extended
by $\le$ and we can also state the above mentioned cancellation law; for the strict relation $<$ both the cancellation
law without limitations and the elementary transfer theorem are valid. Finally for the relation $\le$ we can state a weak
form of
trichotomy: let's write $x \simeq y$ for $x-y \in I_0$ (that is $\st{x}=\st{y}$), then for every $x$, $y \in \ER$
\[
   x \simeq y \quad\text{or}\quad x<y \quad\text{or}\quad y<x.
\]
Anyway neither $\preceq$ nor $\le$ are order relations, as we can see taking $x_t := t \cdot \sin \frac{1}{t}$ which
is not comparable with $y=0$.\\
\indent We conclude this section giving a brief indication of some other possible operations and properties of $\ER$.
First of all we can consider the absolute value: it is a well defined function for which the usual order properties
still hold (use the transfer theorem for inequalities), but for which the following ones are valid too
\begin{align*}
      x \succeq 0 & \quad\Longleftrightarrow\quad  |x| = x \\*[1.5 mm]
      x \preceq 0 & \quad\Longleftrightarrow\quad  |x| = -x \\*[1.5 mm]
      |x| = 0 & \quad\Longleftrightarrow\quad  x = 0.
\end{align*}
Moreover we can consider powers and logarithms of strictly positive (w.r.t. $\le$) extended reals (note that obviously
the square root is not well defined on $D$ therefore the last limitation cannot be eliminate). For these operations
are still valid the usual algebraic and order properties: for example if $y$ is strictly positive and $z > 1$, then we
have
\[
    x \ge y \Longrightarrow \log_z(x) \ge \log_z(y).
\]
\section{The cartesian closure of $\F$}
\label{sec: TheCartesianClosureOfF} In this section we shall define the basic constructions which will lead us to the
notion of $\Cn$ \emph{space} and $\Cn$ \emph{function}. They represent the most general kind of spaces and functions
extendible with our infinitesimal points. Any $\C^n$ manifold is a $\Cn$ space too, and the category $\Cn$ of all
$\Cn$ spaces is cartesian closed, hence it contains several infinite-dimensional spaces, e.g. that formed by all the
usual $\C^n$ functions between two manifolds. It is important to note that, exactly as in \cite{Che, Sou, Fr-Kr,
Kr-Mi, Mo-Re}, the category $\Cn$ contains many ``pathological'' spaces; actually $\Cn$ works as a ``cartesian closed
universe'' and we will see that, like in \cite{Koc, Lav, Mo-Re}, the particular \emph{infinitesimally linear $\Cn$
spaces} have the best properties and will work as a good substitute of manifolds.

The ideas used in this section arise from analogous ideas of \cite{Che} and \cite{Fr-Kr}; actually $\CInfty$ is the
category of diffeological spaces (see \cite{Sou} and references therein).

We present the construction starting from a concrete category $\F$ of topological spaces (which satisfies few
conditions) and embedding it in a cartesian closed category $\bar \F$. We will call $\bar\F$ {\em the cartesian
closure of\/} $\F$. We need this generality because we shall use it to define both domain and codomain of the
extension functor $\ext(-) \;:\; \Cn \freccia \ECn$ starting from two different categories $\F$. The problem to
generalize the definition of $\ER$ to a functor $\ext(-)$ can also be seen from the following point of view: now it is
natural to define a tangent vector as a map
\[
   t : D \freccia \ext{M}.
\]
But we have to note that: $t$ has to be ``regular'' in some sense, hence we need some kind of geometric structure both
on $D$ and $\ext{M}$; the ideal $D$ is not of type $\ext{M}$ for some manifold $M$ because the only standard real
number in $D$ is $0$; the definition of $\ext{M}$ has to generalize $\ER$. We shall define structures on $D$ and
$\ext{M}$
so that $D$, $\ext{M} \in \ECn$, hence we shall define the concept of tangent vector so that $t \in \ECn(D,\ext{M})$.\\
[2mm] \textbf{Hypotheses:}
\begin{enumerate}
    \item $\F$ \emph{is a subcategory of} \textbf{Top} \emph{which contains all the constant maps and all the open subspaces
    $U \subseteq H$ (with
    the induced topology) of every $H \in \F$ with their inclusion $i : U \hookrightarrow H \in \Fs{UH}:=\F(U,H)$.}
\end{enumerate}
\emph{\noindent In the following $|-|\;:\; \F \freccia \Set$ is the forgetful functor which associate to any $H \in F$ its
support set $|H| \in \Set$. Moreover with $\Top{H}$ we will call the topology of $H$ and with $(U \prec H)$ the
subspace of $H$ induced on the open set $U \in \Top{H}$.
\begin{enumerate}
\addtocounter{enumi}{1}
    \item $\F$ is closed with respect to restrictions to open sets, that is if $f \in \Fs{HK}$, $U$ and $V$ are open in
    $H$, $K$ resp. and $f(U)\subseteq V$, then $f|_U \in \F(U\prec H, V \prec K)$;
    \item Every topological space $H \in \F$ has the following ``sheaf property'': let $H$, $K\in \F$, $(H_i)_{i \in I}$
    an open cover of
    $H$ and $f : |H| \freccia |K|$ a map such that $\forall i \in I: f|_{H_i} \in \Fs{H_i K}$, then $f \in \Fs{HK}$.
\end{enumerate}}
For the construction of the domain of the extension functor we want to consider a category $\F$ which permits to embed
 finite dimensional manifolds in $\Cn$. To this aim we will set $\F=\ORn$, the category with objects open sets
 $U \subseteq \R^u$ (with the induced topology), for some $u \in \N$, and with hom-set the usual $\C^n(U,V)$ of $\C^n$ functions between the open
 sets $U \subseteq \R^u$ and $V \subseteq \R^v$. What type of category $\F$ we have to choose depends on the setting we
 need: e.g. in case we want to consider manifolds with boundary we have to take the analogous of the above
 mentioned category $\ORn$ but with objects open set $U \subseteq \R^u_{+}=\{x \in \R^u \,|\, x_u \ge 0\}$.\\
The basic idea to define a $\Cn$ space $X$ (which faithfully generalizes the notion of manifold) is to substitute the
notion of chart with a family of mappings $d : H \freccia X$ with $H \in \F$. E.g. for $\F = \ORn$ these mappings are
of type $d : U \freccia X$ with $U$ open in some $\R^u$, thus they can be thought as $u$-dimensional figures on $X$.
Hence a $\Cn$ space can be thought as a support set and the specification of all the finite-dimensional figures on the
space itself. Generally speaking we can think $\F$ as a category of ``types of figures''. Always considering the case
$\F = \ORn$, we can also think $\F$ as a category which represents ``a well known notion of regular space and regular
function'': with the cartesian closure $\bar \F$ we want extend this notion to a more general type of spaces (e.g.
spaces of mappings). In the diffeological setting \cite{Che, Sou} a figure $d: U \freccia X$ is called a plot on $X$.\\
We are trivially generalizing both the work of \cite{Fr-Kr, Kr-Mi}, where only curves as types of figures are
considered, and the notion of diffeology in which $\F=\ORInfty$. This generalization permit to obtain in an easy way the
cartesian closedness of $\bar\F$, and thus to have at our disposal a general instrument $\F \mapsto \bar\F$ very
useful in the construction e.g. of the codomain of the extension functor $\ext(-)$, where we will choose a different $\F$.
In the sequel we will frequently
use the notation $f \cdot g := g \circ f$ for the composition of maps so as to facilitate the lecture of diagrams, but
we will continue to evaluate functions ``on the right'' hence $(f \cdot g)(x)=g(f(x))$. Objects and arrows of $\bar \F$
generalize the same notions of the diffeological setting.
\begin{definition}
If $X$ is a set, then we say that $(\D,X)$ is an object of $\bar \F$ if $\D = \{\Ds{H}\}_{\sss H \in \F}$ is a family
with
\[
  \Ds{H} \subseteq \Set (|H|,X).
\]
We indicate with the notation $\Fs{JH} \cdot \Ds{H}$ the set of all the compositions $f \cdot d$ of functions $f \in
\Fs{JH}$ and $d \in \Ds{H}$. The family $\D$ has finally to satisfy the following conditions:
\begin{enumerate}
\label{def: BarFobject}
    \item $\Fs{JH} \cdot \Ds{H} \subseteq \Ds{J}$.
    \item $\Ds{H}$ contains all the constant maps $d : |H| \freccia X$.
    \item Let $H \in \F$, $(H_i)_{i \in I}$ an open cover of $H$ and $d : |H| \freccia X$ a map such that
    $d|_{H_i} \in \Ds{(H_i\prec H)}$, then $d \in \Ds{H}$.
\end{enumerate}
Finally we set $|(\D, X)| := X$.
\end{definition}
For the condition \emph{1.} we can think $\Ds{H}$ as the set of all the regular functions defined on the ``well
known'' object $H \in \F$ and with values in the new space $X$; in fact this condition says that the set of figures
$\Ds{H}$ is closed with respect to re-parametrization with $f \in \Fs{JH}$. Condition \emph{2.} is the above mentioned
sheaf property and asserts that to be a figure has a local character depending on $\F$. We will frequently write $d
\In{H} X$ to indicate that $d \in \Ds{H}$ and we can read it ``$d$ is a figure of $X$ of type $H$'' or ``$d$ belong to
$X$ at the level $H$'' or ``$d$ is a generalized element of $X$ of type $H$'' or, finally, ``$(d,U)$ is a plot of
$X$''. This kind of arrows is important to obtain cartesian closure, whereas we shall further use arrows of kind $X
\freccia |H|$ to extend
these spaces with new infinitesimal points.\\
The definition of arrow $f : X \freccia Y$ between two spaces $X$, $Y \in \bar \F$ is the usual one for diffeological
spaces, that is $f : |X| \freccia |Y|$ takes, through composition, generalized elements $d \In{H} X$ of type $H$ in
the domain to generalized elements of the same type in the codomain $f(d):= d \cdot f \In{H} Y$. Note that we have $f:
X \freccia Y$ in $\bar \F$ iff $\forall H \forall x \In{H} X :\; f(x) \In{H}Y$, moreover $X=Y$ iff $\forall H \forall
d:\; d \In{H} X \Leftrightarrow d \In{H} Y$. These and many other properties justify the notation
$\In{H}$ and the name ``generalized element''.\\
With these definitions $\bar \F$ becomes a category. Note that it is, in general, in the second Grothendieck universe
because $\D$ is a family indexed in the set of objects of $\F$ (this is not the case for $\F=\ORn$ which is a set).

The simplest $\bar \F$-object is $\bar K :=(\Fs{(-)K},|K|)$ for $K \in \F$, and for it we have that $d:\bar K \freccia
X$ iff $d \In{K} X$, $\F(H,K)=\bar \F(\bar H,\bar K)$. Therefore $\F$ is fully embedded in $\bar \F$ if $\bar H = \bar
K$ implies $H=K$; e.g. this is true if the given category $\F$ verifies the following hypothesis
\[
  |H|=|K|=S \text{\ \ and\ \ } H \xfrecciad{1_S} K \xfrecciad{1_S} H \then H=K.
\]
E.g. this is true for $\F=\ORn$. Another way to construct an object of $\bar \F$ is to generate it starting from a
given family $\Ds{H}^{\sss{0}} \subseteq \Set(|H|,X)$, for any $H \in \F$, closed with respect to constant functions.
We will indicate this space with $(\F \cdot \D^{\sss{0}}, X)$ and its figures are, locally, compositions $f \cdot d$
with $f \in \Fs{HK}$ and $d \in \Ds{K}^{\sss{0}}$. More precisely $\delta \In{H} (\F \cdot \D^{\sss{0}}, X)$ iff
$\delta : |H| \freccia X$ and for every $h \in |H|$ there exist an open neighborhood $U$ of $h$ in $H$, $K \in \F$, $d
\in \Ds{H}^{\sss{0}}$ and $f : (U \prec H) \freccia K$ in $\F$ such that $\delta|_U = f \cdot d$. \\
On each space $X \in \bar \F$ we can put the final topology $\Top{X}$ for which any figure $d \In{H} X$ is continuous,
that is a subset $U \subseteq |X|$ is in $\Top{X}$ iff $d^{-1}(U) \in \Top{H}$ for any $H \in \F$ and any $d \In{H}
X$. With respect to this topology any arrow of $\bar
\F$ is continuous and we still have the given $\Top{H}$ in the space $\bar H$, that is $\Top{H} = \Top{\bar H}$.\\
Open subsets $U$ on a space $X$ will serve us, e.g., as domains for arrows of type $U \freccia \R^k$. These maps,
which trivially generalize the notion of chart and that we will call ``observables on $X$'', will permit us to define
the extension functor $\ext(-)$.
\subsection{Categorical properties of the cartesian closure}
\label{sec: PropertiesOfCartesianClosure} We shall now examine subobjects in $\bar \F$ and their relation with the
restriction of functions, after which we shall analyse completeness, co-completeness and cartesian closure of $\bar
\F$.
\begin{definition}
Let $X \in \bar\F$ and $S \subseteq |X|$, then we define
\[
    (S \prec X) := (\D,S)                                   \label{def: subspace}
\]
where
\[
    d \in \Ds{H} \DIff d \;:\; |H| \freccia S \e{and} d \cdot i \In{H}X.
\]
Here $i : S \hookrightarrow |X|$ is the inclusion map. We will call $(S \prec X)$ ``the subspace induced on $S$ by
$X$''.
\end{definition}
Using this definition only it's very easy to prove that $(S \prec X) \in \bar\F$ and that its topology contains the
induced topology. Moreover $\Top{(S \prec H)} \subseteq \Top{X}$ if $S$ is open, hence in this case we have on $(S
\prec X)$ the induced topology. Finally we have the following
\begin{theorem}
Let $f : X \freccia Y$ be an arrow of $\bar\F$ and $U$, $V$ subsets of $|X|$ and $|Y|$ respectively, such that $f(U)
\subseteq V$, then
\[
    (U \prec X) \xfrecciad{f|_{\sss{U}}} (V \prec Y) \e{in} \bar\F.
\]
\end{theorem}
Obviously it is easy to state and prove that any $X \in \bar\F$ has the sheaf property. Using our notation for
subobjects we can prove the following useful and natural properties directly from definition \ref{def: subspace}
\begin{itemize}
  \item $(U \prec \bar H)=\overline{(U \prec H)}$ for $U$ open in $H \in \F$
  \item $i : (S \prec X) \hookrightarrow X$ is the lifting of the inclusion $i : S \hookrightarrow |X|$ from $\Set$ to
$\bar \F$
  \item $(|X| \prec X) = X$
  \item $(S \prec (T \prec X)) = (S \prec X)$ \ \ if \ \ $S \subseteq T \subseteq |X|$
  \item $(S \prec X) \times (T \prec Y) = (S \times T \prec X \times Y)$.
\end{itemize}
These properties imply that the relation $X \subseteq Y$ iff $|X| \subseteq |Y|$ and $(|X| \prec Y)=X$ is a partial
order. Note that this relation is stronger to say that the inclusion is an arrow, because it asserts that $X$ and the
inclusion verify the universal property of $(|X| \prec Y)$, that is $X$ is a subobject of $Y$.\\
Completeness and co-completeness are analyzed in the following theorem. For its standard proof see \cite{Fr-Kr} for a
similar theorem.
\begin{theorem}
\label{thm: limit in Cn}
Let $(X_i)_{i \in I}$ be a family of objects in $\bar\F$ and $p_i : |X| \freccia |X_i|$
arrows in $\Set$ $\forall\, i \in I$. Define
\[
    d \In{H} X \DIff d \;:\; |H| \freccia |X| \e{and} \forall\, i \in I \pti d \cdot p_i \In{H} X_i
\]
then $\,\displaystyle(X \xfreccia{p_i} X_i)_{i \in I}\,$ is a lifting of $\,\displaystyle (|X| \xfreccia{p_i}
|X_i|)_{i \in I}\,$ in $\bar\F$.\\
Whereas if $j_i : |X_i| \freccia |X|$ are arrows in  $\Set \forall\, i \in I$ and
\[
    \forall\, x \in |X| \,\, \exists\,\, i \in I \,\, \exists\, x_i \in X_i \pti x=j_i(x_i)
\]
then defining $d \In{H} X$ iff $d: |H| \freccia |X|$ and for every $h \in |H|$ there exist an open neighborhood $U$ of
$h$ in $H$, $i \in I$ and $\delta \In{U} X_i$ s.t. $d|_U = \delta \cdot j_i$, we have that $\,\displaystyle(X_i
\xfreccia{j_i} X)_{i \in I}\,$ is a co-lifting of $\,\displaystyle (|X_i| \xfreccia{j_i} |X|)_{i \in I}\,$ in
$\bar\F$.
\end{theorem}
Directly from the definitions it is easy to prove that on quotient spaces we exactly have the quotient topology and
that on any product we have a topology stronger than the product topology.\\
Finally if we define
\[
  \Ds{H} := \{ d : |H| \freccia \bar\F(X,Y) \;|\; \bar H \times X \xfrecciad{d^\vee} Y \e{in} \bar\F \}
\]
(we are using the notations $d^\vee(h,x) := d(h)(x)$ and $\mu^\wedge(x)(y) := \mu(x,y)$) then $\langle \D, \bar\F(X,Y)
\rangle =: Y^X$ is an object of $\bar \F$. With this definition, see e.g. \cite{Che} or \cite{Fr-Kr}, it is easy to
prove that $\bar \F$ is cartesian closed, i.e. that the $\bar \F$-isomorphism $(-)^\vee$ realizes
\[
   (Y^X)^Z \simeq Y^{Z \times X}.
\]
\section{The category $\Cn$}
\subsection{Observables on $\Cn$ spaces and separated spaces}
The most natural way to apply the results of previous section for our aims is to set $\F=\Man$, that is to consider
directly the cartesian closure of the category of finite dimensional $\C^n$ manifolds (we shall not formally assume
any hypothesis on the topology of a manifold because we will never need it in the following; moreover if not
differently specified, with the word ``manifold'' we will always mean ``finite dimensional manifold''). We shall not
follows this idea for several reasons; we will set instead $\Cn := \overline{\ORn}$, that is the cartesian closure of
the category $\ORn$ of open sets and $\C^n$ arrows. For $n=\infty$ this gives exactly diffeological spaces \cite{Che,
Sou}. As we noted before $\overline{\Man}$ is in the second Grothendieck universe and, essentially for simplicity,
from this point of view the choice $\F = \ORn$ is better. In spite of this choice it is natural to expect, and in fact
we will prove it, that the category of finite-dimensional manifolds is faithfully contained in $\Cn$. Another reason
for our definition of $\Cn$ is that in this way the category of $\Cn$ spaces and arrows is more natural to accept and
to work in with respect to $\overline{\Man}$; hence ones again a reason of simplicity. We will see that manifolds
modelled in convenient vector spaces (see \cite{Kr-Mi}) are faithfully embedded in $\Cn$, hence our choice to take
finite dimensional objects in the definition of $\Cn$ is not restrictive from this
point of view.\\
Now we pay attention to another type of maps which go ``on the opposite direction'' with respect to figures $d : K
\freccia X$. As mentioned above we shall use them to introduce new infinitesimal points for any $X \in \Cn$.
\begin{definition}
Let $X$ be a $\Cn$ space, then we say that
\[
    UK \text{\emph{ is a zone (in $X$)}}
\]
iff $U \in \Top{X}$ is open in $X$ and $K \in \ORn$. Moreover we say that
\[
    c \text{\emph{ is an observable on }$UK$} \qquad \text{or} \qquad c \InUp{UK} X
\]
iff $c \;:\; (U \prec X) \freccia \bar K$ in $\Cn$.
\end{definition}
Remember that for any open set $K$ the $\Cn$ space $\bar K$ is
\[
    \bar K = (\C^n(-,K), K),
\]
hence composition of figures $d \In{H}X$ with observables $c \InUp{UK}$ gives ordinary $\C^n$ maps: $d|_S \cdot c \in
\C^n(S,K)$, where $S:=d^{-1}(U)$.\\
From our previous theorems it follows that $\Cn$ functions $f : X \freccia Y$ take observables on the codomain to
observables on the domain i.e.:
\begin{equation}
    c \InUp{UK} Y \then f|_{\sss{S}} \cdot c \InUp{SK} X,                           \label{thm: fPreservesObs}
\end{equation}
where $S := f^{-1}(U)$. Therefore isomorphic $\Cn$ spaces have isomorphic sets of figures and observables.\\
Generalizing through the observables the equivalence relation \ref{def: sim} to generic $\Cn$ spaces, we will have to
study the following condition, which is connected with the faithfulness of the extension itself.
\begin{definition}
\label{def: xIdentifiedy}
If $X \in \Cn$ and $x$, $y \in |X|$, then we write
\[
    x \asymp y
\]
and we read it ``$x$\emph{ and }$y$\emph{ are identified in }$X$'', iff for every zone $UK$ and every $c \InUp{UK} X$ we
have
\begin{enumerate}
    \item $x \in U \iff y \in U$
    \item $x \in U \then c(x) = c(y).$
\end{enumerate}
Moreover we say that $X$\emph{ is separated} iff $x \asymp y$ implies $x=y$ for any $x$, $y \in |X|$.
\end{definition}
Observe that if two points are identified in $X$ then a generic open set contains the first iff it contains the second
too (take a constant observable), and from (\ref{thm: fPreservesObs}) that $\Cn$ functions $f : X \freccia Y$ preserve the
relation $\asymp$:
\[
    x \asymp y \text{ in } X \then f(x) \asymp f(y) \text{ in } Y.
\]
Trivial examples of separated spaces can be obtained considering the objects $\bar U$ with $U \in \ORn$, or taking
subobjects of separated spaces. But the full subcategory of separated $\Cn$ spaces has other good enough properties.
\begin{theorem}
The category of separated $\Cn$ spaces is complete and admits co-products. Moreover if $X$, $Y$ are separated then
$Y^X$ is separated too, and hence separated spaces form a cartesian closed category.
\end{theorem}
{\em \ \ Sketch of the proof:\/} We only do some considerations about co-product, because it is easy to prove that
products and equalizers of separated spaces are separated too. Let us consider a family $(\mathcal{X}_i)_{i \in I}$ of
separated spaces with support sets $X_i := |\mathcal{X}_i|$. Constructing their sum in $\Set$
\[
    X := \sum_{i \in I} X_i
\]
\[
    j_i \;:\; x \in X_i \longmapsto (x,i) \in X,
\]
from the completeness of $\Cn$ we can lift it to a co-product $\displaystyle(\mathcal{X}_i \xfreccia{j_i}
\mathcal{X})_{i \in I}$. To prove that $\mathcal{X}$ is separated we take two points $x$, $y \in X=|\mathcal{X}|$
identified in $\mathcal{X}$. These points are of the form $x=(x_r,r)$ and $y=(y_s,s)$, with $x_r \in X_r, y_s \in X_s$
and $r,s \in I$. We want to prove that $r$ and $s$ are necessarily equal. In fact from the definition of figures of
$\mathcal{X}$ (Theorem \ref{thm: limit in Cn}) we have that
\[
    A \in \Top{\mathcal{X}} \quad \iff \quad \forall \, i \in I \pti j_i^{-1}(A) \in \Top{\mathcal{X}_i},
\]
and hence $X_r \times \{r\}$ is open in $\mathcal{X}$ and $x \asymp y$ implies
\[
    (x_r,r) \in X_r \times \{r\} \iff (y_s,s) \in X_r \times \{r\} \ee{hence} r=s.
\]
Hence $x=y$ iff $x_r$ and $y_s=y_r$ are identified in $\mathcal{X}_r$ and this is a consequence of the following
facts:
\begin{enumerate}
\item if $U$ is open in $\mathcal{X}_r$ then $U \times \{r\}$ is open in $\mathcal{X}$; \item if $c \InUp{UK}
\mathcal{X}_r$, then $\gamma(x,r) := c(x)\,\, \forall\, x \in U$ is an observable of $\mathcal{X}$ on $U \times
\{r\}$.
\end{enumerate}
Now let us consider exponential objects. If $f$, $g \in |Y^X|$ are identified, to prove that they are equal is
equivalent to prove that $f(x)$ and $g(x)$ are identified in $Y$ for any $x$. To obtain this conclusion is sufficient
to consider that the evaluation in $x$ i.e. $\eps_x : \varphi \in |Y^X| \longmapsto \varphi(x) \in |Y|$ is a $\Cn$ map
and hence from any observable $c \InUp{UK}Y$ we can always obtain the observable $\eps_x|_{\sss U^{\prime}} \cdot c
\InUp{U^{\prime}K} Y^X$ where $U^{\prime} := \eps_x^{-1}(U)$.\qed Finally let's consider two $\Cn$ spaces such that
the topology $\Top{X\times Y}$ is equal to the product of the topologies $\Top{X}$ and $\Top{Y}$. Then if $x$,
$x^\prime \in |X|$ and $y$, $y^\prime \in |Y|$ it is easy to prove that we have $x \asymp x'$ in $X$ and $y
\asymp y'$ in $Y$ iff $(x,y) \asymp (x', y')$ in $X \times Y$.
\subsection{Manifolds as objects of $\Cn$}
We can associate in a very natural way a $\Cn$ space $\bar M$ to any manifold $M \in \Man$ with the following
\begin{definition}
    Define $|\bar M|:=|M|$ and for every $H \in \ORn$
    \[
        d \In{H} \bar M \DIff d \in \Man(H,M).
    \]
\end{definition}
We obtain a $\Cn$ space with the same topology of the starting manifold. Moreover the observables of $\bar M$ are the
most natural that one could expect, in fact it is very easy to prove that
\[
    c \InUp{UK} \bar M \quad\iff\quad c \in \Man(U,K).
\]
Hence it is clear that $\bar M$ is separated, because charts are observables of the space. The following theorem says
that the passage from $\Man$ to $\Cn$ that we are considering is a full embedding and therefore it says that $\Cn$ is
a non-trivial generalization of the notion of manifold which include infinite-dimensional spaces too.
\begin{theorem}
\label{thm: ManifoldsEmbedded}
Let $M$ and $N$ be $\C^n$ manifolds, then
\begin{enumerate}
    \item $\bar M = \bar N \then M=N$
    \item $\displaystyle\bar M \xfrecciad{f} \bar N \e{in} \Cn \ \iff\ M \xfrecciad{f} N \e{in} \Man$.
\end{enumerate}
Hence $\Man$ is fully embedded in $\Cn$.
\end{theorem}
\emph{Proof of 1)}: If $(U,\varphi)$ is a chart on $M$, then $\varphi^{-1}|_A : A:=\varphi(U) \freccia M$ is a figure
of $\bar M$, that is $\varphi^{-1}|_A \In{A} \bar M = \bar N$. But if $\psi : U \freccia \psi(U) \subseteq \R^k$ is
a chart of $N$, then it is also an observable of $\bar N$, and composition of figures and observables gives ordinary
$\C^n$ maps, that is the atlases of $M$ and $N$ are compatible.\\
\emph{Proof of 2)}: We use the same ideas as above and moreover that $\varphi^{-1}|_A \In{A} \bar M$ implies
$\varphi^{-1}|_A \cdot f \In{A} \bar N$. Finally we can compose this $A$-figure of $\bar N$ with a chart (observable)
of $N$ obtaining an ordinary $\C^n$ map. \qed Directly from the definitions we can prove that for two manifolds we
also have
\[
    \overline{M \times N} = \bar M \times \bar N.
\]
This property is useful to prove the affirmations done in the following examples.
\addtocontents{toc}{\protect\newpage}
\subsection{Examples}
\begin{enumerate}
\item Let $M$ be a $\C^\infty$ manifold modelled on convenient vector spaces (see \cite{Kr-Mi}). We can define $\bar M$
analogously as above, saying that $d \In{H} \bar M$ iff $d : H \freccia M$ is a smooth map between $H$ (open in some
$\R^h$) and $M$. In this way smooth curves on $M$ are exactly the figures $c \in_{\R} \bar M$ of type $\R$ in $\bar
M$. On $M$ we obviously think the so called natural topology, that is the identification topology with respect to some
smooth atlas, which is also the final topology with respect to all smooth curves and hence is also the final topology
$\Top{\bar M}$ with respect to all figures of $\bar M$. More easily with respect to the previous case of finite
dimensional manifolds, it is possible to study observables, obtaining that $c \InUp{UK} \bar M$ iff $c : U \freccia K$
is smooth as a map between manifolds modelled on convenient vector spaces. Moreover if $(U,\varphi)$ is a chart of $M$
on the convenient vector space $E$, then $\varphi : (U \prec \bar M) \freccia (\varphi(U) \prec \bar E)$ is $\CInfty$.
Using these results it is easy to prove the analogous of Theorem \ref{thm: ManifoldsEmbedded} for the category of
manifolds modelled in convenient vector spaces. Hence also classical smooth manifolds modelled in Banach spaces are
embedded in $\CInfty$.
\item It is possible to prove that the following applications, frequently used e.g. in calculus of variations, are smooth,
that is they are arrows of $\CInfty$.
\begin{enumerate}
\item The operator of derivation:
\begin{align*}
  \partial_i :\ \C^\infty(\R^n,\R^k) &\freccia \C^\infty(\R^n,\R^k) \\
              u      &\longmapsto \frac{\partial u}{\partial x_i}
\end{align*}
\item The integral operator:
\begin{align*}
  i :\ \C^\infty(\R^2,\R) &\freccia \C^\infty(\R,\R)\\
     u &\longmapsto \displaystyle\int_a^b u(-,s) \diff{s}
\end{align*}
\item Using the previous examples we can prove that the classical operator of calculus of variations is smooth
\[
  \mathcal{I}(u)(t) := \int_a^b F[u(t,s),\partial_2 u(t,s),s] \diff{s}
\]
\[
  \mathcal{I} : \C^\infty(\R^2,\R^k) \freccia \C^\infty(\R,\R),
\]
where the function $F : \R^k \times \R^k \times \R \freccia \R$ is smooth.
\end{enumerate}
\item Because of cartesian closedness set-theoretical operations like the following are examples of $\Cn$ arrows:
    \begin{itemize}
        \item composition:
        \[
            (f,g) \in B^A \times C^B \;\; \mapsto \;\; g \circ f \in C^A
        \]
        \item evaluation:
        \[
            (f,x) \in Y^X \times X \;\; \mapsto \;\; f(x) \in Y
        \]
        \item insertion:
        \[
            x \in X \;\; \mapsto \;\; (x,-) \in (X \times Y)^Y
        \]
    \end{itemize}
\item Inversion between smooth manifolds modelled on Banach spaces
    \[
        (-)^{-1} : f \in \text{Diff}(N,M) \;\; \mapsto \;\; f^{-1} \in \text{Diff}(M,N)
    \]
is a smooth mapping, where $\text{Diff}(M,N)$ is the subspace of $N^M$ given by the diffeomorphisms between $M$ and
$N$. So $(\text{Diff}(M,M),\circ)$ is a (generalized) Lie group. To prove that $(-)^{-1}$ is smooth let's consider a
figure $d \In{U}\text{Diff}(N,M)$, then $f:=(d \cdot i)^\vee : U \times N \freccia M$, where $i : \text{Diff}(N,M)
\hookrightarrow M^N$ is the inclusion, is an ordinary smooth function between Banach manifolds. We have to prove that
$g:=[d \cdot (-)^{-1} \cdot j]^\vee : U \times M \freccia N$ is smooth, where $j : \text{Diff}(M,N) \hookrightarrow
N^M$. But $f[u,g(u,m)]=m$ and $\text{\bf D}_2f(u,n)=\text{\bf D}[d(u)](n)$ hence the conclusion follows from the
implicit function theorem because $d(u) \in \text{Diff}(N,M)$.
\item Since the category $\Cn$ is complete, we can also
have $\Cn$ spaces with singular points like e.g. the equalizer $\{ x \in X \,|\, f(x)=g(x) \}$. Any algebraic curve is
in this way a $\CInfty$ separated space too. \item Another type of space with singular points is the following. Let
$\varphi \in \C^n(\R^k,\R^m)$ and consider the subspace $([0,1]^k \prec \R^k)$, then $(\varphi([0,1]^k) \prec \R^m)
\in \Cn$ is the deformation in $\R^m$ of the hypercube $[0,1]^k$.
\item Let $C$ be a continuum body, $I$ the interval
for time, and $\cal E$ the 3-di\-men\-sional Euclidean space. We can define on $C$ a natural structure of
$\boldsymbol{\C}^\infty$ space. For any point $p \in C$ let $p_r(t) \in \cal E$ be the position of $p$ at time $t$ in
the frame of reference $r$; we define figures of type $U$ on $C$ ($U \in \ORn$) the functions $d : U \freccia C$ for
which the following application
\begin{align*}
  \tilde{d} : U \times I &\freccia \mathcal{E} \\
             (u,t) &\longmapsto d(u)_r(t)
\end{align*}
is smooth. For example if $U=\R$ then we can think $d : \R \freccia C$ as a curve traced on the body and parameterized
by $u \in \R$. Hence we are requiring that the position $d(u)(t)$ of the particle $d(u) \in C$ varies smoothly with the
parameter $u$ and the time $t$. This is a generalization of the continuity of motion of any point of the body (take $d$
constant). This smooth (that is diffeological) space will be separated, as an object of $\boldsymbol{\C}^\infty$, if
different points of the body cannot have the same motion:
\[
   p_r(-)=q_r(-) \then p=q  \qquad \forall p,q \in C.
\]
The configuration space of $C$ can be viewed (see \cite{Wan}) as the space
\[
     M := \sum_{t \in I} M_t \ee{with} M_t \subseteq \mathcal{E}^C
\]
and so, for the categorical properties of $\boldsymbol{\C}^\infty$ the spaces $\mathcal{E}^C$, $M_t$ and $M$ are
always objects of $\boldsymbol{\C}^\infty$ as well. With this structure the motion of $C$:
\begin{align*}
  \mu :\ &C \times I \freccia \mathcal{E}  \\
       {}&\h{.055}(p,t) \longmapsto p(t)
\end{align*}
is a smooth map. Note that to obtain these results we need neither $M_t$ nor $C$ be manifolds, but only the possibility
to associate to any point $p$ of $C$ a motion $p_r(-) : I \freccia \mathcal{E}$. If we had the possibility to develop
differential geometry for these spaces too we would have the possibility to obtain many results of continuum mechanics
for body which cannot be naturally represented using a manifolds or with infinite-dimensional configuration space.
Moreover in the next section we will see how to extend any $\boldsymbol{\C}^\infty$ space with
infinitesimal points, so that we can also consider infinitesimal sub-bodies of $C$.
\end{enumerate}
\section{The extension of $\Cn$ spaces and functions}
Now we want to extend any $\Cn$ space and any $\Cn$ function by means of our ``infinitesimal points''. First of all we
will have to extend to a generic space $X$ the notion of nilpotent path. Afterward we shall use the observables to
generalize the equivalence relation $\sim$ (see Definition \ref{def: sim}) using the following idea
\[
    \varphi(x_t) = \varphi(y_t) + \text{o}(t) \ee{with} \varphi \InUp{UK}X.
\]
In this point the main problem is to understand how to relate $x$, $y$ with the domain $U$ of $\varphi$. In the
subsequent sections we will also prove some results that will conduct us toward the theorem $\ext{(M \times N)} \simeq
\ext{M} \times \ext{N}$ with $M,N$ manifolds. The fact that this useful theorem is not proved for generic $\Cn$ spaces
is due to the fact that the topology on a product between $\Cn$ spaces is generally stronger than the product
topology.
\subsection{Nilpotent paths}
If $X$ is a $\Cn$ space, then using $\Top{X}$ we can define the set $\C_0(X)$ of all the maps $x : \R \freccia X$
continuous at the origin $t=0$. Because any $\Cn$ function $f$ is continuous we have $f \circ x \in \C_0(Y)$ if
$x \in \C_0(X)$.\\
If $U$ is open in $X$ then on the subspace $(U \prec X)$ we have the induced topology and from this it follows that
\begin{theorem}
    \label{thm: LimitWithObs}
    Let $X$ be a $\Cn$ space and $x \in \C_0(X)$. Take an observable $\varphi \InUp{UK}X$ with $x(0) \in U$, then
    \[
        \lim_{t \to 0} \varphi(x_t) = \varphi(x_0).
    \]
\end{theorem}
As many other concepts we will introduce, the notion of nilpotent map is defined by means of observables.
\begin{definition}
    Let $X$ be a $\Cn$ space and $x \in \C_0(X)$, then we say that $x$ \emph{is nilpotent (rel. $X$)} iff for every zone
    $UK$ of $X$ and every obsevable $\varphi \InUp{UK}X$ we have
    \[
        x(0) \in U \then \exists k \in \N : \| \varphi(x_t) - \varphi(x_0) \|^k = \text{\rm o}(t).
    \]
    Moreover
    \[
        \No{X}:=\No{}(X):= \{ x \in \C_0(X)\,\, |\,\, x \text{\rm \ is nilpotent}\}.
    \]
\end{definition}
Because of property (\ref{thm: fPreservesObs}), if $f \in \Cn(X,Y)$ and $x \in \No{X}$ then $f \circ x \in \No{Y}$, that
is $\Cn$ functions preserve nilpotent maps. In case of a manifold $M$, a map $x : \R \freccia |M|$ is nilpotent iff we
can find a chart $(U,\varphi)$ on $x_0$ such that $\| \varphi(x_t) - \varphi(x_0) \|^k = \text{\rm o}(t)$ for some $k
\in \N$.\\
Finally we enunciate the relations between product manifolds and nilpotent paths. For the (standard) proof is
essential to observe that $\Top{\bar M \times \bar N} = \Top{\overline{M \times N}} = \Top{M \times N}$ and thus on
the product $\bar M \times \bar N$ of $\Cn$ spaces we exactly have the product topology.
\begin{theorem}
    \label{thm: NilpotentAndProduct}
    Let $M,N$ be manifolds and $x : \R \freccia |M|$, $y : \R \freccia |N|$, then
    \[
        x \in \No{\bar M} \e{and} y \in \No{\bar N} \iff
        (x,y) \in \No{\bar M \times \bar N}.
    \]
    Here $(x , y)_t := (x_t,y_t)$.
\end{theorem}
\subsection{The extension of spaces and functions}
\label{sec: extensions}
\begin{definition}
    Let $X$ be a $\Cn$ space and $x,y \in \No{X}$ then we say that
    \[
        x \sim y \e{in} X \ee{or simply} x=y \e{in} \ext{X}
    \]
    iff for every zone $UK$ of $X$ and every observable $\varphi \InUp{UK} X$ we have\\*[4 mm]
    1)\ \ $x_0 \in U \iff y_0 \in U$ \\*[2 mm]
    2)\ \ $x_0 \in U \then \varphi(x_t) = \varphi(y_t) + \text{\rm o}(t).$
\end{definition}
Obviously we write $\ext{X} := \No{X}/\!\sim$ and $\ext{f}(x) := f \circ x$ if $f \in \Cn(X,Y)$ and $x \in \ext{X}$.
We prove the correctness of the definition of $\ext{f}$ in the following:
\begin{theorem}
    If $f \in \Cn(X,Y)$ and $x=y$ in $\ext{X}$ then $\ext{f}(x)=\ext{f}(y)$ in $\ext{Y}$.
\end{theorem}
\Proof Take a zone $VK$ in $Y$ and an observable $\psi \InUp{VK} Y$, then from continuity of $f$, $U := f^{-1}(V) \in
\Top{X}$. We can thus apply hypothesis $x=y$ in  $\ext{X}$ with the zone $UK$ and the observable $\varphi := f|_U
\cdot \psi \InUp{UK} X$. From this it follows the conclusion noting that $f \circ x, f \circ y \in \No{Y}$ and $x_0
\in U$ iff $f(x_0) \in V$. \qed Using Theorem \ref{thm: LimitWithObs} we can note that $x=y$ in $\ext{X}$ implies that
$x_0$ and $y_0$ are identified in $X$ (Definition \ref{def: xIdentifiedy}) and thus using constant maps $\hat x (t) :=
x$ we obtain an injection $\hat{\scriptstyle(-)} : |X| \freccia \ext{X}$ if the space $X$ is separated. Therefore if
$Y$ is separated too, $\ext{f}$ is really an extension of $f$. Finally note that $\ext{(-)}$ preserves
compositions and identities.\\
{\em Example:} If $X=M$ is a $\C^n$ manifold then we have $x \sim y$ in $M$ iff there exists a chart $(U,\varphi)$ of
$M$ such that
\begin{enumerate}
    \item \ $x_0$, $y_0 \in U$
    \item \ $\varphi(x_t) = \varphi(y_t) + \text{\rm o}(t)$.
\end{enumerate}
Moreover the previous conditions do not depend on the chart $(U,\varphi)$. In par\-ticu\-lar if $X=U$ is an open set
in $\R^k$, then $x \sim y$ in $U$ is simply equivalent to the limit relation $x(t) = y(t) + \text{o}(t)$; hence if $i
: U \hookrightarrow \R^k$ is the inclusion map, it's easy to prove that $\ext{i} : \ext{U} \freccia \ER^k$ is
injective. As in Theorem \ref{thm: DandOpen} we will always identify $\ext{U}$ with $\ext{i}(\ext{U})$, so we simply
write $\ext{U} \subseteq \ER^k$. Using this equivalent way to express the relation $\sim$ on manifolds, we can see
that $(x,y) = (x',y')$ in $\ext{(M \times N)}$ iff $x=x'$ in $\ext{M}$ and $y = y'$ in $\ext{N}$. From this conclusion
and from Theorem \ref{thm: NilpotentAndProduct} we can prove that the following applications
\[
    \alpha_{\sss MN}:= \alpha \,\,:\,\, ([x]_{\sim},[y]_{\sim}) \in \ext{M} \times \ext{N} \longmapsto
    [(x,y)]_{\sim} \in \ext{(M \times N)}
\]
\[
    \beta_{\sss MN}:= \beta \,\,:\,\, [z]_{\sim} \in \ext{(M \times N)} \longmapsto
    ([z \cdot p_M]_{\sim}, [z \cdot p_N]_{\sim}) \in \ext{M} \times
    \ext{N}
\]
(for clarity we have used the notation with the equivalence classes) are well-defined bijections with $\alpha^{-1} =
\beta$ (obviously $p_{\sss M}, p_{\sss N}$ are the projections). We will use the first one of them in the following
section with the temporary notation $\langle p,x \rangle :=\alpha(p,x)$, hence $f \langle p,x \rangle =
f(\alpha(p,x))$ for $f : \ext{(M \times N)} \freccia Y$. This simplifies our notations but it permits to avoid the
identification of $\ext{M} \times \ext{N}$ with $\ext{(M \times N)}$ until we will have proved that $\alpha$ and
$\beta$ are arrows of the category $\ECn$.
\section{The category of extended spaces and the extension functor}
\subsection{Motivations:}
Up to now every $\ext{X}$ is a simple set only. Now we want to use the general passage from a category $\F$ to its
cartesian closure $\bar \F$ so as to put on any $\ext{X}$ some kind of useful structure. Our aim is to obtain in this
way a new cartesian closed category $\bar \F =: \ECn$ and a functor ``extension'' $\ext{(-)} : \Cn \freccia \ECn$.
Therefore we have to choose $\F$, that is what will be the types of figures of $\ext{X}$. It may seem very natural to
take $\ext{g} : \ext{U} \freccia \ext{V}$ as arrow in $\F$ if $g : U \freccia V$ is in $\ORn$ (in \cite{Gio} we
followed this way). The first problem in this idea is that, e.g.
\[
    \ER \xfrecciad{\ext{f}} \ER \then \ext{f}(0)=f(0) \in \R,
\]
hence there cannot exist a constant function of type $\ext{f}$ to a non-standard value and so we cannot satisfy the
closure of $\F$ with respect to generic constant functions (see hypotheses on $\F$ in Section \ref{sec:
TheCartesianClosureOfF}). But we can make further considerations about this problem so as to motivate better the
choice of $\F$. The first one is that we surely want to have the possibility to lift, using cartesian closedness, maps
simple as the sum between extended reals:
\[
    s : (p,q) \in \ER \times \ER \freccia p+q \in \ER.
\]
Hence $s^\wedge(p): q \in \ER \freccia p+q \in \ER$ must be an arrow of $\ECn$. Note that it is nor constant neither
of type $\ext{f}$ because $s^\wedge(p)(0)=p$ and $p$ could be an extended real.\\
The second consideration is about $\alpha$: if we want to have $\alpha$ as an arrow of $\ECn$, then in the following
situation we have to obtain a $\ECn$ arrow again
\begin{gather*}
    \ER \times \ER \xfrecciad{p \times 1_{\ER}} \ER \times \ER \xfrecciad{\alpha} \ext{(\R \times \R)} \xfrecciad{\ext{g}}
    \ER\\
    (t,s) \longmapsto (p,s) \longmapsto \langle p,s \rangle \longmapsto \ext{g}\langle p,s \rangle
\end{gather*}
(where $p \in \ER$ and $g \in \C^n(\R^2,\R)$). The idea we shall follow is exactly to take as arrows of $\F$ maps that
locally are of type $\delta(s)=\ext{g}\langle p,s \rangle$, where $p$ works as a parameter of $\ext{g}\langle -,-
\rangle$. Obviously in this way $\delta$ could also be a constant map to an extended value (take as $g$ a projection).
Frequently one can find maps of type $\ext{g}\langle p,- \rangle$ in informal calculations in physics or geometry.
Actually they simply are $\C^n$ maps with some fixed parameter $p$, which could be an infinitesimal distance (e.g. in
the potential of the electric dipole, see below), an infinitesimal coefficient associated to a metric (like in the
formula given at the beginning), or
considering a side of an infinitesimal surface.\\
Note the importance of $\alpha$ to perform passages like the following
\begin{align*}
    {}& M \times N \xfrecciad{f} Y&{}                            & \e{in $\Cn$} \\
    {}& \ext{(M \times N)} \xfrecciad{\ext{f}} \ext{Y} &{}       & \e{in $\ECn$} \\
    {}& \ext{M} \times \ext{N} \xfrecciad{\ext{f}} \ext{Y} &{}   & \e{in $\ECn$ (identification via $\alpha$)}\\
    {}& \ext{M} \xfrecciad{\ext{f}^\wedge} \ext{Y}^{\ext{N}} &{} & \e{using cartesian closedness.}
\end{align*}
This motivates the choice of arrows in $\F$, but there is a second problem about the choice of its objects. Take a
manifold $M$ and an arrow $t : D \freccia \ext{M}$ in $\ECn$. Whatever this will mean we want to think $t$ as a
tangent vector applied either to a standard point $t(0) \in M$, and in this case it is a standard tangent vectors, or
to an extended one, $t(0) \in \ext{M} \setminus M$. Roughly speaking this is the case if we can write
$t(h)=\ext{g}\langle p,h \rangle$ for some $g$, $p$. If we want to obtain this equality it is useful to have both
$1_D$ as a figure of $D$
\[
   1_D \In{D} D \then t \In{D}\ext{M},
\]
and maps of type $\ext{g} \langle p,-\rangle : D \freccia \ext{M}$ as figures of $\ext{M}$. Therefore it would be useful
to have $D$ as an object of $\F$. But $D$ is not the extension of a standard subset of $\R$, thus what will be the objects
 of $\F$? We will take generic subsets $S$ of $\ext{(\R^s)}$ with the topology $\Top{S}$ generated by $\mathcal{U} =
 \ext{U} \cap S$ for $U$ open in $\R^s$ (in this case we will say that the open set $\mathcal{U}$ \emph{is defined by}
 $U$ \emph{in} $S$). These are the motivations to introduce $\F$ by means of the following
\begin{definition}
    We call $\SERn$ the category whose objects are subsets $S \subseteq \ext{(\R^s)}$, for some $s$ which depends on $S$,
    and with the previous topology $\Top{S}$.
    If $S \subseteq \ext{(\R^s)}$ and $T \subseteq \ext{(\R^t)}$ then we say that
    \[
        S \xfrecciad{f} T \e{in} \SERn
    \]
    iff $f$ maps $S$ in $T$ and for every $s \in S$ we can write
    \[
        f(x) = \ext{g}\langle p,x \rangle \quad \forall x \in \mathcal{V}
    \]
    for some
    \begin{align*}
        {}& p \in  \ext{(\R^\text{\sf p})} \\
        {}& \mathcal{U} \text{ open neighborhood of $p$ defined by $U$ in } \ext{(\R^\text{\sf p})}\\
        {}& \mathcal{V} \text{ open neighborhood of $s$ defined by $V$ in } S\\
        {}& g \in \C^n(U \times V,\R^t).
    \end{align*}
    Moreover we will consider on $\SERn$ the forgetful functor given by the inclusion
    $|-| : \SERn \hookrightarrow \text{\bf Set}$.
\end{definition}
It is easy to prove that $\SERn$ and the functor $|-|$ verify the hypotheses on $\F$ (see Section \ref{sec:
TheCartesianClosureOfF}), hence we can define
\[
    \ECn := \overline{\SERn}.
\]
Each object of $\ECn$ is called an ``extended ($\Cn$) space''.
\subsection{The extension functor}
Now the problem is: what extended spaces could we associate to sets like $\ext{X}$ or $D$? For any subset $Z \subseteq
\ext{X}$ we call $\ext{(ZX)}$ the extended space generated on $Z$ (see Section \ref{sec: TheCartesianClosureOfF}) by
the following set of figures $d : T \freccia Z$ (where $T \subseteq \ext{(\R^t)}$)
\begin{equation}
\label{def: GeneratingSetForExtension}
    \begin{aligned}
        d \in \Ds{T}^0(Z) \DIff & \text{$d$ is constant or we can write} \\
        {}& \text{$d = \ext{h}|_T$ for some $h \In{V}X$ such that $T \subseteq \ext{V}$.}
    \end{aligned}
\end{equation}
Thus in the non-trivial case we start from a standard figure $h \In{V}X$, where the extension of $V$ contains $T$; we
extend this figure obtaining $\ext{h} : \ext{V} \freccia \ext{X}$, and finally the restriction $\ext{h}|_{\sss T}$ is
a generating figure if it maps $T$ in $Z$.\\
Using this definition we call (with some abuses of language)
\begin{gather*}
    \ext{X}:=\ext{(\ext{X}X)}\\
    D :=\ext{(D\R)}\\
    \ER:=\ext{(\ext{\R}\R)}\\
    \ER^k:=\ext{(\ext{(\R^k)}\R^k)}\\
    D_k := \ext{(D_k\R^k)}.
\end{gather*}
We will call $\ext{(ZX)}$ \emph{the extended space induced on $Z$ by $X$}. We can now study the extension functor:
\begin{theorem}
    Let $f \in \Cn(X,Y)$ and $Z$ a subset of $\ext{X}$ with $\ext{f}(Z) \subseteq W \subseteq \ext{Y}$, then in $\ECn$
    we have that
    \[
        \ext{(ZX)} \xfrecciad{\ext{f}|_Z} \ext{(WY)}.
    \]
    Therefore $\ext{(-)} : \Cn \freccia \ECn$.
\end{theorem}
\Proof Take a figure $\delta \In{S} \ext{(ZX)}$ in the domain. We have to prove that $\delta \cdot \ext{f}|_Z$ locally
factors through $\SERn$ and $\D^0(W)$. Hence taking $s \in S$ we can write $\delta|_U = f_1 \cdot d$ where $U$ is an
open neighborhood of $s$, $f_1 \in \SERn(U \prec S,T)$ and $d \in \Ds{T}^0(Z)$. We omit the trivial case $d$ constant,
hence we can suppose to have, using the same notation as above, $d=\ext{h}|_T: T \freccia Z$ with $h \In{V} X$.
Therefore
\[
   ( \delta \cdot \ext{f}|_Z )|_U = f_1 \cdot \ext{h}|_T \cdot \ext{f}|_Z = f_1 \cdot \ext{(h f)}|_T.
\]
But $h f \In{V} Y$ and so $( \delta \cdot \ext{f}|_Z )|_U = f_1 \cdot d_1$, where $d_1 := \ext{(h f)}|_T \in
\Ds{T}^0(W)$, which is the conclusion.\qedNoNewLine
\subsection{The isomorphisms $\alpha$ and $\beta$}
We want to prove that the above mentioned bijective applications $\alpha$ and $\beta$ are arrows of $\ECn$. To
simplify the proof we will use the following preliminary results. The first one is a general property of the extension
$\F \mapsto \bar{\F}$.
\begin{lemma}
\label{lem: ProductOfGeneratedObjects}
    Suppose that $\F$ admits finite products, and for every objects $K$, $J$ an isomorphism
    $$
    \begin{CD}
        \beta := \beta_{\sss{KJ}} : \overline{K \times J} @>>\sim> \bar{K} \times \bar{J} \e{in} \bar{\F}.
    \end{CD}
    $$
    Now let $Z$, $X$, $Y \in \bar{\F}$ with $X$ and $Y$ generated by $\D^{\sss X}$ and $\D^{\sss Y}$ respectively. Then
    we have
    \[
        X \times Y \xfrecciad{f} Z \e{in} \bar{\F}
    \]
    iff for any $K$, $J \in \F$ and $d \in \Ds{K}^{\sss X}$, $\delta \in \Ds{J}^{\sss Y}$ we have
    \[
        \beta \cdot (d \times \delta) \cdot f \In{K \times J} Z
    \]
\end{lemma}
The second Lemma asserts that $\F=\SERn$ verifies the hypotheses of the previous one.
\begin{lemma}
\label{lem: ProductInSERn}
    The category $\SERn$ admits finite products and the above mentioned isomorphisms $\beta_{\sss KJ}$. Moreover let
    $M$, $N$ be $\C^n$
    manifolds, and $h \In{V}M$, $l \In{V'}N$ with $K \subseteq \ext{V}$ and $J \subseteq \ext{V'}$, then
    \[
        \beta_{\sss KJ} \cdot ( \ext{h}|_{\sss K} \times \ext{l}|_{\sss J} ) \cdot \alpha_{\sss MN} =
        \ext{( h \times l )}|_{\sss K \times J}
    \]
\end{lemma}
The proofs are a direct effect of the given definitions.
\begin{theorem}
    Let $M$, $N$ be $\C^n$ manifolds, then in $\ECn$ we have
    \[
        \ext{(M \times N)} \simeq \ext{M} \times \ext{N}.
    \]
\end{theorem}
\Proof Note that in the statement each manifold is identified with the corresponding $\Cn$ space $\bar{M}$. Hence we
mean $\ext{M}=\ext{\bar{M}}=\ext{(\ext{M}\bar{M})}$. In proving that $\alpha$ is a $\ECn$ arrow we can use the Lemma
\ref{lem: ProductOfGeneratedObjects} because of Lemma \ref{lem: ProductInSERn} and considering that $\ext{M}$ and
$\ext{N}$ are generated by $\D^0(\ext{M})$ and $\D^0(\ext{N})$. Because these generating sets are defined using a
disjunction (see definition \ref{def: GeneratingSetForExtension}) we have to check four cases depending on $d$ and
$\delta$. In the first one we have $d=\ext{h}|_{\sss K} \in \Ds{K}^0(\ext{M})$ and $\delta=\ext{l}|_{\sss J} \in
\Ds{J}^0(\ext{N})$ (we are using the same notations of the previous Lemma \ref{lem: ProductInSERn}). Thus
\[
    \beta_{\sss KJ} \cdot (d \times \delta) \cdot \alpha =
    \beta_{\sss KJ} \cdot (\ext{h}|_{\sss K} \times \ext{l}|_{\sss J}) \cdot \alpha =
    \ext{(h \times l)}|_{\sss K \times J}.
\]
That is $\beta_{\sss KJ} \cdot (d \times \delta) \cdot \alpha$ is a generating element in $\ext{(M \times N)}$, and so
it is also a figure. In the second case let's suppose $\delta$ constant to $n \in \ext{N}$, take a chart $l^{-1} : U
\freccia \R^\text{\sf p}$ on $\st{n}=n_0 \in N$ and let $p:=\ext{l^{-1}}(n)$, $W:=\ER^\text{\sf p}$. Then for any $k
\in K$ and $j \in J$ we can write
\begin{align}
\label{eq: AlphaBetaDelta}
    \alpha\{(d \times \delta)[ \beta_{\sss{KJ}}(\langle k,j \rangle)]\}&=\alpha[\ext{h}(k), \ext{l}(p)] \notag \\
   {}&= \{ \beta_{\sss{K W}}\cdot [\ext{h}|_{\sss K} \times \ext{l}|_{\sss J}] \cdot \alpha \}\langle k,p \rangle\notag\\
   {}&= \ext{(h \times l)}|_{\sss K \times W} \langle k,p \rangle
\end{align}
where we have used once again the equality of Lemma \ref{lem: ProductInSERn}. Thus let's call $\tau : \langle k,j
\rangle \in |K \times J| \mapsto \langle k,p \rangle \in |K \times W|$, so that we can write (\ref{eq:
AlphaBetaDelta}) as
\[
   \beta_{\sss{KJ}} \cdot (d \times \delta) \cdot \alpha = \tau \cdot \ext{(h \times l)}|_{\sss K \times W}
\]
But $\ext{(h \times l)}|_{\sss K \times W}$ is a generating figure of $\ext{(M \times N)}$ and $\tau$ is an arrow of
$\SERn$, and this proves that $\beta_{\sss{KJ}} \cdot (d \times \delta) \cdot \alpha \In{K \times J} \ext{(M \times
N)}$. The remaining cases are either trivial or analogous to the last one.\\
For $\beta_{\sss MN}$ the proof is simpler, and it suffices to note that e.g. $\alpha \cdot \ext{p}_{\sss M}$ is the
projection on $\ext{M}$; hence the conclusion follows from the fact that $\ext{p}_{\sss M}$ and $\ext{p}_{\sss N}$ are
arrows of $\ECn$. \qed
In the following we shall always use $\alpha$ to identify these type of spaces $\ext{M}\times\ext{N} = \ext{(M \times
N)}$.
\subsection{Figures of extended spaces}
\label{sec: FiguresOfExtendedSpaces}
In this section we want to understand better the figures of the extended
space $\ext{(ZX)}$; we will use these results later, for example when we will study the embedding of $\Man$ in $\ECn$.\\
From the general definition of $\bar\F$--space generated by $\D^0$ given in section \ref{sec: TheCartesianClosureOfF},
a figure $\delta \In{S} \ext{(ZX)}$ can be locally factored as $\delta|_\mathcal{V} = f d$ through an arrow $f$ of
$\SERn$ and a generating function $d \in \Ds{T}^0(Z)$; here $\mathcal{V}=\mathcal{V}(s)$ is an open neighborhood of
the fixed $s \in S$. Hence, either $\delta|_\mathcal{V}$ is constant (if $d$ is constant) or we can write
$d=\ext{h}|_T$ and $f=\ext{g}(p,-)$ so that
\[
   \delta(x)=d[f(x)]=\ext{h}[\ext{g}(p,x)]=\ext{(g h)}(p,x) \quad \forall x \in \mathcal{B},
\]
where $\mathcal{B}=\ext{B} \cap \mathcal{V}$ and $A \times B$ is an open neighborhood of $(\st{p},\st{s})$. Therefore
we can write
\[
   \delta(x) = \ext{\gamma}(p,x) \quad \forall x \in \mathcal{B}
\]
with $\gamma := g|_{A \times B} \cdot h \in \Cn(A \times B,X)$. Thus figures of $\ext{(ZX)}$ are locally necessarily
either constant map or a natural generalization of the maps of $\SERn$, that is ``parameterized extended $\Cn$
arrows''. Using the properties of $\ECn$ and of its arrow $\alpha_{\R^\text{\sf p} \R^s}$ it is easy to prove that
these conditions are sufficient too. Moreover if $X=M$ is a manifold, the condition ``$\delta|_\mathcal{V}$ constant''
can be omitted. In fact if $\delta|_\mathcal{V}$ is constant to $m \in Z \subseteq \ext{M}$, then taking a chart
$\varphi$ on $\st{m} \in M$ we can write $\delta(x)=m=\ext{\gamma}(p,x)$, where $p=\ext{\varphi}(m)$ and
$\gamma(x,y)= \varphi^{-1}(x)$. Using these results we can see that $\ext{(ZX)}=(Z \prec \ext{X})$. Hence if we take
$Z \subseteq \ER^z$, we have three coincident ways to see it as an extended space: $\bar{Z}=\ext{(Z\R^z)}=(Z \prec
\ER^z)$ (here $\bar{Z}$ is the general passage from an object $H \in \F$ to $\bar{H} \in \bar{\F}$). E.g. if $f :
\ER^z \freccia \ext{X}$ is a $\ECn$ arrow, then we also have $f : \overline{\ER^z} \freccia \ext{X}$ and so $f
\In{\ER^z} \ext{X}$ and locally we can write $f$ either as a constant function or, with the usual notations, as
$f(x)=\ext{\gamma}(p,x)$. For functions $f : I \freccia \ext{X}$ defined on some set $I \subseteq I_0$ of
infinitesimals which contains $0 \in I$, these two alternatives are globally true instead of locally only.\\
We close this section enunciating the following properties of the extension functor:
\begin{enumerate}
    \item If $X \subseteq Y$ in $\Cn$ (see section \ref{sec: PropertiesOfCartesianClosure}) and $|X|$ is open in $Y$,
    then $\ext{X} \subseteq \ext{Y}$ in $\ECn$ and $\ext{X}$ is open in $\ext{Y}$.
    \item In the same hypotheses as above, if $Z \subseteq |\ext{X}|$ then $(Z \prec \ext{X}) = (Z \prec \ext{Y})$.
    \item Let $f : X \freccia Y$ and $Z \subseteq Y$ in $\Cn$, with $|Z|$ open in $Y$. Moreover define
    $\ext{f}^{-1}(\ext{Z}):= (\ext{f}^{-1}(|\ext{Z}|) \prec \ext{X})$ and
    $f^{-1}(Z):= (f^{-1}(|Z|) \prec X)$. Then $\ext{[f^{-1}(Z)]}=\ext{f}^{-1}(\ext{Z})$ as extended spaces.
\end{enumerate}
\subsection{The embedding of manifolds in $\ECn$}
\label{sec: EmbeddingInExtendedSpaces}
If we consider a $\Cn$ space $X$, we have just seen that we have the possibility to associate an extended space to any
subset $Z \subseteq \ext{X}$. Thus if $X$ is separated we can put a structure of $\ECn$ space on the set $|X|$ of ordinary
points of $X$, by means of $\bar{X} := \ext{(|X|X)}=(|X| \prec \ext{X})$. Intuitively $X$ and $\bar X$ seem very similar,
and in fact we have
\begin{theorem}
    \label{thm: EmbeddingOfSepInExtendedCn}
    Let $X$, $Y$ be $\Cn$ separated spaces, then
    \begin{enumerate}
        \item \ $\bar{X} = \bar{Y} \then X=Y$
        \item \ $\bar{X} \xfrecciad{f} \bar{Y} \e{in} \ECn \ \iff\ X \xfrecciad{f} Y \e{in} \Cn.$
    \end{enumerate}
    Hence $\Cn$ separated spaces are fully embedded in $\ECn$, and so is $\Man$.
\end{theorem}
\Proof 1) The equality $\bar{X} = \bar{Y}$ immediately implies the equality of support sets $|X|=|Y|$. We consider now
a generalized element $d\In{H} X$ where $H$ is an open set of $\R^h$. Taking the extension of $d$ and then the
restriction to ordinary points only we obtain
\begin{equation}
    (H \prec \ext{\bar{H}}) \xfrecciad{\ext{d}|_{\sss H}} (|X| \prec \ext{X})=\bar{X}=\bar{Y}.
\end{equation}
But $(H \prec \ext{\bar{H}})=(H \prec \ER^h)=\ext{(H \R^h)}=\bar{H}$, hence
\[
   \ext{d}|_{\sss H} = d : \bar{H} \freccia \bar{Y} \e{in} \ECn
\]
and so $d \In{H} \bar{Y}$. Therefore for every $s \in H$ either $d$ is constant in some open neighborhood
$\mathcal{V}$ of $s$ defined by $V$, or, using the usual notations, we can write
\[
    d(x) = \ext{\gamma}(p,x) \text{ in } \ext{Y}\quad \forall x \in \mathcal{V}=\ext{V} \cap H = V \cap H.
\]
Hence for every $x \in V \cap H$ we have $\st{d(x)} \asymp \st{[\gamma(p,x)]}$ in $Y$, and so we can write $d(x) =
\gamma(p_0,x)$ because $Y$ is separated and $x \in V \cap H \subseteq \R^h$ is standard. Therefore $d|_{V \cap H}$ is
a $Y$-valued arrow of $\Cn$ defined in a neighborhood of the fixed $s$. The conclusion thus follows from the sheaf
property of
$Y$.\\
2 $\Rightarrow$) From the proof of 1) we have seen that if $d \In{H} X$ then $d \In{H} \bar{X}$. Hence $f(d) \In{H}
\bar{Y}$. But once again from the passages of 1) we have seen that this implies that $f(d) \In{H} Y$.\\
2 $\Leftarrow$) It is sufficient to extend $f$, to restrict it to standard points only, and finally to consider
that our spaces are separated. \qed
An immediate corollary of this theorem is that the extension functor is another full embedding for separated spaces.
\begin{corollary}
    Let $X,Y$ be $\Cn$ separated spaces, then
    \begin{enumerate}
        \item \ $\ext{X} = \ext{Y}  \then X=Y$
        \item \ If\ \ $ \ext{X} \xfrecciad{f} \ext{Y}$ in $\ECn$ and $f(|X|) \subseteq |Y|$ then
              \[
                    X \xfrecciad{f|_{|X|}} Y \e{in} \Cn
              \]
        \item \ $\ext{X} \xfrecciad{\ext{f}} \ext{Y} \e{in} \ECn\ \iff\ X \xfrecciad{f} Y \e{in} \Cn$
        \item \ If $f$, $g : X \freccia Y$ are $\Cn$ functions, then
              \[
                  \ext{f}=\ext{g} \then f=g.
              \]
    \end{enumerate}
\end{corollary}
{\em Proof of 1):\/} We have to prove that the support sets of $X$ and $Y$ are equal, but this is trivial if we
take standard parts
\[
    \{ \st{x} \,\,|\,\, x \in \ext{X} \} = |X| = \{ \st{x} \,\,|\,\, x \in \ext{Y} \} = |Y|.
\]
Hence $\bar{X}=(|X| \prec \ext{X}) = (|Y| \prec \ext{Y})=\bar{Y}$. The other properties stated in the corollary are
proved considering the previously seen passages which uses the restriction to ordinary points, in case preceded by an
application of the extension functor. \qedNoNewLine
\subsection{The generalized derivation formula in $\ECn$}
\label{sec: GeneralizedDerivationFormula} In this section we want to explore the possibility to use the derivation
formula through the use of observables $\varphi \InUp{UK} X$. Precisely we start from a generic $\ECn$ function $f : D
\freccia \ext{X}$ with $f(0) \in \ext{U} := \ext{(U \prec X)}$ (not the extension of a classical one, that is $f$
generally is not of the form $f = \ext{g}|_D$) and we study the validity of the formula for the function
$\ext{\varphi}(f(-))$ defined in $D$ and with values in the $\ER$ module $\ER^k$. Note that the result is not trivial
just because the function $f$ generally is not of the form $\ext{g}|_D$ but of a more general type. First of all we
prove that the previous composition is well defined, that is the following generalization of Theorem \ref{thm:
DandOpen}
\begin{theorem}
\label{thm: DandOpenGeneralized}
    Let $X$ be a $\Cn$ space and $U \in \Top{X}$ an open set. Let $f : D \freccia \ext{X}$ be a $\ECn$ function with
    $f(0) \in \ext{U}$, then $f(h) \in \ext{U}$ for every $h \in D$.
\end{theorem}
\Proof From the hypothesis on $f$ it follows that $f \In{D} \ext{X}$ because $D=\overline{D}$. Hence, considering
that $0 \in D \subset I_0$ and the results of section \ref{sec: FiguresOfExtendedSpaces}, we can globally say that either
$f$ is constant, and the proof is trivial, or we can write the equality $f(h) = \ext{\gamma}(p,h)$ in $\ext{X}$ for
every $h \in D$. For the sake of clarity let $y:=f(h)$, thus taking standard parts (that is evaluating at $t=0$)
\begin{equation}
\label{eq: StOf_y_and_StOff0}
    \st{y} \asymp \st{[\ext{\gamma}(p,h)]}=\gamma(p_0,0)=\st{[\ext{\gamma}(p,0)]} \asymp \st{f(0)}.
\end{equation}
But $f(0) \in \ext{U}$, from which $\st{f(0)} \in U$ and so $\st{y} \in U$ from the previous relation (\ref{eq:
StOf_y_and_StOff0}). Hence $y_t \in U$ for $t$ small and moreover $y \in \No{U}$ because $y=f(h) \in \No{X}$ and
because on $U=(U \prec X)$ we have the induced topology.\qedNoNewLine
\begin{theorem}
    \label{thm: GeneralizedDerivationFormula}
    Let $X$ be a $\Cn$ space and $\varphi \InUp{U\R^k} X$ an observable. Let $f$ be as above, then there exists one and
    only one pair
    \[
        a \in \ER^k \ee{and} b \in \R^k
    \]
    such that
    \[
        \forall h \in D \pti \ext{\varphi}(f(h)) = a + h \cdot b.
    \]
\end{theorem}
\Proof Omitting as usual the trivial case in which $f$ is constant, from the previous proof we have seen that we can
write
\[
    \forall h \in D \pti f(h) = \ext{\gamma}(p,h) \e{in} \ext{X}.
\]
Therefore from the definition of equality in $\ext{X}$
\[
    \varphi[f(h)_t] = \varphi[\gamma(p_t,h_t)] + t \cdot \sigma_1(t) \quad \forall^0 t,
\]
with $\lim_{t \to 0} \sigma_1(t)=0$. But the function $\psi := \varphi[\gamma(-,-)]$ is an ordinary $\C^n$ function,
hence we can use the Taylor formula  ($n$ is at least $1$) to obtain
\begin{gather*}
    \psi(p_t,h_t) = \psi(p_t,0) + h_t \cdot \partial_2\psi(p_t,0) + \sigma_2(t) \\
    \lim_{\substack{t \to 0\\ h_t \ne 0}} \frac{\sigma_2(t)}{h_t} = 0.
\end{gather*}
Hence if we define
\begin{align*}
    {}& a := \ext{\varphi}[\ext{\gamma}(p,0)] \in \ER^k \\
    {}& b := \partial_2 \psi(p_0,0) \in \R^k
\end{align*}
then substituting
\begin{multline*}
    \varphi[f(h)_t] - a_t - h_t \cdot b = \\
    = h_t \cdot [\partial_2 \psi(p_t,0) - \partial_2 \psi(p_0,0)] + \sigma_2(t) + t \cdot \sigma_1(t) = \text{o}(t)
\end{multline*}
This proves that $\ext{\varphi}[f(h)]=a+h \cdot b$. To prove uniqueness of $a$ is sufficient to set $h=0$; for $b$ is
sufficient to note that if
\[
    \forall h \in D \pti h \cdot b = h \cdot \beta
\]
then setting $h_t = t$, from the equality in $\ER^k$ we quickly obtain $b=\beta$. \qed Using the generalized derivation
formula we can extend Theorem \ref{thm: DerivationFormula} to non-standard points $x \in \ER$. It suffices to
consider the function $f(x + \,\boldsymbol{\cdot}) : D \freccia \ER$
which is an arrow in $\ECn$ and to which we can apply Theorem \ref{thm: GeneralizedDerivationFormula}.\\
We will denote with $\varphi'(f)$ the unique $b$ in the previous theorem so that we can formulate the following
result, in which is stated that the generalized derivation formula determine uniquely the function $f$.
\begin{theorem}
    \label{thm: UniquenessDerivationFormula}
    Let $X \in \Cn$ and $f,g : D \freccia \ext{X}$ in $\ECn$, with $f(0)=g(0)$. Moreover we     assume that
    $\varphi'(f) = \varphi'(g)$ for every $\varphi \InUp{U\R^k} X$ with $f(0) \in    \ext{U}$, then f=g.
\end{theorem}
\Proof Fix an $h \in D$ and for simplicity let $y := f(h)$ and $z := g(h)$. From the proof of a previous theorem we
have seen that $\st{y}\asymp\st{f(0)}$, but $f(0)=g(0)$ in $\ext{X}$ hence $y_0 \asymp \st{f(0)} \asymp \st{g(0)}
\asymp z_0$. Now we consider an observable $\varphi \InUp{U\R^k}X$ and from $y_0 \asymp z_0$ we deduce that
\[
    z_0 \in U \iff y_0 \in U.
\]
Hence if we assume that $y_0 \in U$ then $f(0)_t \in U \ \forall^0 t$ and $f(0) \in \ext{U}$. Thus from the hypotheses
of the theorem it follows that $\varphi'(f) = \varphi'(g)$ and hence the generalized derivation formula
implies that $\ext{\varphi}(f(h))=\ext{\varphi}(g(h))$, that is $\varphi(y_t)=\varphi(z_t) + {\rm o}(t)$. \qed
In the case that $X$ is a manifold, to have the equality $f=g$ is sufficient to find a chart $(U,\varphi)$ with $f(0) \in
\ext{U}$ and for which $\varphi'(f) = \varphi'(g)$ (see the example in Section \ref{sec: extensions}).
\section{Examples}
We started this article defining in a very simple way an extension $\ER$ of the real field containing nilpotent
infinitesimals. Afterwards, generalizing diffeological spaces, we introduced a cartesian closed embedding $\Cn$ of
$\Man$ to which we generalized the definition of $\ER$ obtaining the category $\ECn$. The aim of this article is to
introduce the foundations of this theory of infinitesimals, leaving its full development in Differential Geometry for
future works. To perform this aim it is important to note the deep analogy between our construction and Synthetic
Differential Geometry (see \cite{Lav} and references therein): frequently we only have to trivially generalize this
work.

The elementary examples listed in the following want to show in a few rows the simplicity of the analytic/algebraic
calculus using nilpotent elements. Here ``simplicity'' means that the dialectic with informal calculations is really
faithful; this is important for future developments
 both as a proof of the flexibility of the new language  and also for researches in artificial intelligence like
automatic differentiation theories. Last but not least it may also be important for didactical or historical researches.
\begin{enumerate}
    \item \textbf{Commutation of differentiation and integration.} Suppose we want \emph{to discover} the derivative of the
    function
    \[
        g(x) := \int_{\alpha(x)}^{\beta(x)} f(x,t) \diff{t} \qquad \forall x \in \R
    \]
    where $\alpha$, $\beta$ and $f$ are $\C^1$ functions. We can see $g$ as a composition of locally lipschitzian
    functions hence we can apply the derivation formula:
    \begin{align*}
      g(x+h) = &\int_{\alpha(x)+h \alpha'(x)}^{\alpha(x)} f(x,t) \diff{t} + h \cdot \int_{\alpha(x)+h
                    \alpha'(x)}^{\alpha(x)}\frac{\partial f}{\partial x}(x,t) \diff{t} + \\
             {}&+ \int_{\alpha(x)}^{\beta(x)} f(x,t) \diff{t} + h \cdot \int_{\alpha(x)}^{\beta(x)}
                        \frac{\partial f}{\partial x}(x,t) \diff{t} + \\
             {}&+ \int_{\beta(x)}^{\beta(x)+ h \beta'(x)} f(x,t) \diff{t} + h \cdot \int^{\beta(x)+h \beta'(x)}_{\beta(x)}
                \frac{\partial f}{\partial x}(x,t) \diff{t}.  \\
    \end{align*}
    Now we use $h^2 = 0$ to obtain e.g.
    \[
        h \cdot \int_{\alpha(x)+h \alpha'(x)}^{\alpha(x)}\frac{\partial f}{\partial x}(x,t) \diff{t} =
        -h^2 \cdot \alpha'(x) \cdot \frac{\partial f}{\partial x}(x,t)=0
    \]
    and
    \[
        \int_{\alpha(x)+h \alpha'(x)}^{\alpha(x)} f(x,t) \diff{t} = -h \cdot \alpha'(x) \cdot f(\alpha(x),t).
    \]
    Treating in an analogous way similar terms we finally obtain the conclusion.
    Note that the final formula comes out by itself so that we have \emph{``discovered''} it and not simply we have proved it.
    \item \textbf{Circle of curvature.} A simple application of the infinitesimal Taylor formula is the parametric
    equation for the circle of curvature, that is the circle with second order osculation with a curve
    $\gamma : [0,1] \freccia \R^3$. In fact if $r \in (0,1)$ and $\dot\gamma_r$ is a unit vector, from the second order
    formula we have
    \begin{equation}
    \label{eq:_CircleOfCurvature}
        \forall h \in D_2 \pti \gamma(r+h) = \gamma_r + h \, \dot\gamma_r + \frac{h^2}{2} \,\ddot \gamma_r =
        \gamma_r+h \,\vec t_r + \frac{h^2}{2}c_r \, \vec n_r
    \end{equation}
    where $\vec n$ is the unit normal vector, $\vec t$ is the tangent one and $c_r$ the curvature. But once again from
    Taylor formula we have $\sin(c h) =  c h$ and $\cos(c h) = 1 - \frac{c^2 h^2}{2}.$ Now it
    suffices to substitute $h$ and $\frac{h^2}{2}$ from these formulas into (\ref{eq:_CircleOfCurvature}) to obtain the
    conclusion
    \[
        \forall h \in D_2 \pti \gamma(r+h) = \left( \gamma_r + \frac{\vec n_r}{c_r} \right) +
        \frac{1}{c_r} \cdot \left[ \sin(c_r h) \vec t_r - \cos(c_r h) \vec n_r \right].
    \]
    In a similar way we can prove that any $f \in \C^\infty(\R,\R)$ can be written $\forall h \in D_k$ as
    \[
        f(h)=\sum_{n=0}^k a_n \cdot \cos(n h) + \sum_{n=0}^k b_n \cdot \sin(n h).
    \]
    \item \textbf{Schwarz's theorem.} Using nilpotent infinitesimals a simple and meaningful proof of Schwarz's theorem
    can be obtained. This simple example aims to show how to manage some differences between our setting and Synthetic
    Differential Geometry (see \cite{Koc, Lav, Mo-Re}).
    Let $f : V \freccia E$ be a $\C^2$ function between Banach spaces and $a \in V$, we want to prove
    that ${\rm d}^2{f}(a) : V \times V \freccia E$ is symmetric. Take
    \begin{align*}
        {}& k \in D_2 \\
        {}& h,j \text{ infinitesimals}\\
        {}& j k h \in D_{\ne 0}
    \end{align*}
    (e.g. we can take $k_t=|t|^{\frac{1}{2}}, h_t=j_t=|t|^{\frac{1}{4}}$). Using $k \in D_2$ we have
    \begin{equation}
        \label{eq: x+h+k}
        \begin{split}
            j \cdot f(x & + h u + k v) =\\
            &= j \cdot \left[ f(x+ h u) + k \, \partial_v f(x + h u) + \frac{k^2}{2}\partial^2_v f(x + h u) \right]\\
            &= j \cdot f(x + h u) + j k \cdot \partial_v f(x + h u)
        \end{split}
    \end{equation}
    where we used the fact that $k^2 \in D$ and $j$ infinitesimal imply $j k^2 = 0$. Now we consider that $j k h \in D$
    hence using Theorem \ref{thm: SecondDerivationFormula} we obtain
    \begin{equation}
        \label{eq: jkh}
        j k \cdot \partial_v f(x + h u) = j k \cdot \partial_v f(x) + j k h \cdot \partial_u(\partial_v f)(x).
    \end{equation}
    But $k \in D_2$ and $j k^2=0$ hence
    \[
        j \cdot f(x + k v) - j \cdot f(x) = j k \cdot \partial_v f(x).
    \]
    Substituting in (\ref{eq: jkh}) and (\ref{eq: x+h+k}) we obtain
    \begin{equation}
        \label{eq: SecondOrderIncrementalRatio}
        \begin{split}
        & j \cdot \left[ f(x + h u + k v) - f(x + h u) - f(x + k v) +f(x) \right] =\\
        & =j k h \cdot \partial_u(\partial_v f)(x).
        \end{split}
    \end{equation}
    The left side of this equality is symmetric in $u,v$, hence changing them we have
    \[
        j k h \cdot \partial_u(\partial_v f)(x) = j k h \cdot \partial_v(\partial_u f)(x)
    \]
    and hence  the conclusion because $j k h \ne 0$ and $\partial_u(\partial_v f)(x),\partial_v(\partial_u f)(x) \in E$.\\
    From (\ref{eq: SecondOrderIncrementalRatio}) it follows directly the classical limit relation
    \[
        \lim_{t \to 0}\frac{f(x + h_t u + k_t v) - f(x + h_t u) - f(x + k_t v) + f(x)}{h_t k_t} = \partial_u(\partial_v
        f)(x).
    \]
    \item \textbf{Electric dipole}. From a Physical point of view an electric dipole is usually defined as
    ``\emph{a pair of charges with opposite sign placed at a distance $d$ very less than the distance  $r$ from the
    observer''}.\\
    Conditions like $r \gg d$ are frequently used in Physic and very often we obtain a correct formalization if we
    ask $d \in \ER$ infinitesimal but $r \in \R \setminus \{0\}$ i.e. $r$ finite. Thus we can define an electric dipole as a pair
    $(p_1,p_2)$ of electric particles, with charges of equal intensity but
    with opposite sign such that their mutual distance at every time $t$ is a first order infinitesimal:
    \begin{equation}
        \label{eq:_DefinitionOfDipole}
        \forall t \pti \vert p_1(t) - p_2(t) \vert =: \vert \vec{d}_t \vert =: d_t \in D.
    \end{equation}
    In this way we can calculate the potential in the point $x$ using the properties of $D$ and using the hypothesis that
    $r$ is finite and not zero. In fact we have
    \[
        \varphi(x) = \frac{q}{4\pi\epsilon_0} \cdot
        \left( \frac{1}{r_1} - \frac{1}{r_2} \right)
        \qquad \qquad
        \vec{r_i} := x - p_i
    \]
    and if $\vec r := \vec{r}_2 - \frac{\vec d}{2}$ then
    \[
        \frac{1}{r_2} = \left( r^2 + \frac{d^2}{4} + \vec r \boldsymbol{\cdot} \vec d \right)^{-1/2}
        = r^{-1} \cdot \left( 1 + \frac{\vec r \boldsymbol{\cdot} \vec d}{r^2} \right)^{-1/2}
    \]
    because for (\ref{eq:_DefinitionOfDipole}) $d^2=0$. For our hypotheses on $d$ and $r$ we have that
    $\displaystyle\frac{\vec r \boldsymbol{\cdot} \vec d}{r^2} \in D$ hence from the derivation formula
    \[
        \left( 1 + \frac{\vec r \boldsymbol{\cdot} \vec d}{r^2} \right)^{-1/2} =
        1 - \frac{\vec r \boldsymbol{\cdot} \vec d}{2 r^2}
    \]
    In the same way we can proceed for $1/r_1$, hence:
    \[
        \varphi(x) = \frac{q}{4\pi\epsilon_0} \cdot \frac{1}{r} \cdot
        \left( 1 + \frac{\vec r \boldsymbol{\cdot} \vec d}{2 r^2} - 1 +
        \frac{\vec r \boldsymbol{\cdot} \vec d}{2 r^2} \right) = \ldots
    \]
    The property $d^2 = 0$ is also used in the calculus of the electric field and for the moment of momentum.
    \item \textbf{Newtonian limit in Relativity.} Another example in which we can formalize a condition like $ r \gg d$
    using the previous ideas is the Newtonian limit in Relativity; in it we can suppose to have
    \begin{itemize}
        \item $\,\forall t \pti v_t \in D_2 \e{and} c \in \R$
        \item $\, \forall x \in M_4 \pti g_{i j}(x)=\eta_{i j}+h_{i j}(x) \e{with} h_{i j}(x) \in D.$
    \end{itemize}
    where $\left( \eta_{i j} \right)_{i j}$ is the matrix of the Minkowski's  metric. This conditions can be interpreted
    as $v_t \ll c$ and $h_{ij}(x) \ll 1$ (low speed with respect to the speed of light and weak gravitational field). In
    this way we have, e.g. the equalities:
    \[
        \frac{1}{\sqrt{\displaystyle 1 - \frac{v^2}{c^2}}} = 1 + \frac{v^2}{2 c^2}
        \ee{and} \sqrt{1 - h_{44}(x)} = 1 - \frac{1}{2}\,h_{44}(x).
    \]
    \item \textbf{Linear differential equations.}  Let
        \begin{gather*}
            L(y):=A_{\sss 0}\frac{\diff{}^{\sss N} y}{\diff{}t^{\sss N}} +
            \ldots + A_{\sss N-1} \frac{\diff{} y}{\diff{}t} + A_{\sss N} \cdot y=0
        \end{gather*}
    be a linear differential equation with constant coefficients. Once again we want \emph{to discover} independent
    solutions in case the characteristic polynomial has multiple roots e.g.
    \[
        (r - r_{\sss 1})^2 \cdot (r - r_{\sss 3}) \cdot \ldots \cdot (r- r_{\sss N})=0.
    \]
    The idea is that in $\ER$ we have $(r - r_1)^2=0$ also if $r = r_1 + h$ with $h \in D$. Thus
    $y(t)={\rm e}^{(r_1 + h)t}$ is a solution too. But ${\rm e}^{(r_1 + h)t}={\rm e}^{r_1 t} + h t \cdot {\rm e}^{r_1 t}$,
    hence
    \begin{align*}
    L\left[{\rm e}^{(r_1 + h)t}\right]&=0\\
                                    {}&=L\left[{\rm e}^{r_1 t} + h t \cdot {\rm e}^{r_1 t} \right]\\
                                    {}&=L\left[{\rm e}^{r_1 t} \right] + h \cdot L\left[t \cdot {\rm e}^{r_1 t} \right]
    \end{align*}
    We obtain $L\left[t \cdot {\rm e}^{r_1 t} \right]=0$, that is $y_1(t)=t \cdot {\rm e}^{r_1 t}$ must be a solution.
    Using $k$-th order infinitesimals we can deal with other multiple roots in a similar way.
\end{enumerate}
\section{Tangent vectors, vector fields and infinitesimally linear spaces}
The use of nilpotent infinitesimals permits to develop many concepts of Differential Geometry in an intrinsic way
without being forced to use coordinates as we shall see in some examples below. In this way the use of charts
becomes specific of stated areas.\\
We can call this kind of intrinsic geometry \emph{Infinitesimal Differential Geometry}.\\
The possibility to avoid coordinates using infinitesimal neighborhood instead permits to perform some generalizations
to more abstract spaces, like spaces of mappings. Even if the categories $\Cn$ and $\ECn$ are very big and not very
much can be said about generic objects, in this section we shall see that the best properties can be formulated for a
restricted class of extended spaces, the infinitesimally linear ones, to which spaces of mappings between manifolds
belong to.

We start from the fundamental idea of tangent vector. It is now natural to define a tangent vector to a space $X \in
\ECn$ as an arrow (in $\ECn$) of type $t : D \freccia X$. Therefore ${\rm T}X := X^D$ and ${\rm T}f(t) := \diff{f}(t)
:= f \circ t$ with projection $\pi : t \in {\rm T}X \mapsto t(0) \in X$ is the tangent bundle of $X$.  Note that using
the absolute value it is also possible to consider ``boundary tangent vectors'' taking $ |D|:= \{ \; |h| :\, h \in
D\}$ instead of $D$, for example at the initial point of a curve or at a side of a closed set. In the following $M \in
\text{\bf Man}^{\infty}=:\ManInfty$ will always be a finite dimensional smooth manifold and we will use the
notation ${\rm T}M$ for ${\rm T}(\ext{M})$.\\
It is important to note that with this definition of tangent vector we obtain a generalization of the classical
notion. In fact $t(0) \in \ext{M}$ and hence the tangent vector $t$ can be applied to an extended point. If we want to
study classical tangent vectors only we have to consider the following $\CInfty$ object
\begin{definition}
\label{def: StndTangentFunctor}
We call ${\rm T}_{\rm st}M$ the $\CInfty$ object with support set
\[
    \left| {\rm T}_{\rm st}M \right| := \{\ext{f}|_D \,:\, f \in \CInfty(\R,M) \},
\]
and with generalized elements of type $U$ (open in $\R^u$)
\[
    d \In{U} {\rm T}_{\rm st}M \DIff d: U \freccia \left| {\rm T}_{\rm st}M \right|
          \e{and} d \cdot i \In{\bar{U}} {\rm T}M,
\]
where $i : \left| {\rm T}_{\rm st}M \right| \hookrightarrow {\rm T}M$ is the inclusion.
\end{definition}
That is in ${\rm T}_{\rm st}M$ we consider only tangent vectors $t=\ext{f}|_D$ obtained as extension of ordinary
smooth functions $f: \R \freccia M$, and we take as generalized elements, functions $d$ with values in ${\rm T}_{\rm st}M$
which in $\ECInfty$ verify $d^{\vee} : \bar{U} \times D \freccia \ext{M}$ (here $U \in {\rm \bf S}\ER^\infty \mapsto
\bar{U} \in \CInfty$ is the general passage from an object $H \in \F$ to $\bar{H} \in \bar{\F}$). Note that, intuitively
speaking, $d$ takes standard element $u \in U \subseteq \R^k$ to standard element $d(u) \in {\rm T}_{\rm st}M$.
\begin{theorem}
\label{thm: StandardTangentVector}
    Let $t \in {\rm T}M$ be a tangent vector, then
\[
    t \in {\rm T}_{\rm st}M \iff t(0) \in M.
\]
\end{theorem}
\Proof If $t=\ext{f}|_D$ then $t(0)=f(0) \in M$. Vice versa if $t(0) \in M$ then take a chart $(U, \varphi)$ on $t(0)$
and apply the generalized derivation formula (Theorem \ref{thm: GeneralizedDerivationFormula}) obtaining
$\ext{\varphi}(t(h))=a + h \cdot b$ for any $h \in D$ and with $a \in \ER^k$, $b \in \R^k$. But
$\ext{\varphi}(t(0))=\varphi(t(0))=a$ because $t(0) \in M$. Hence $a \in \R^k$ is standard and we can write
$t(h)=\ext{\varphi}^{-1}(a + h \cdot b)=: \ext{f}|_D(h)$. \qed In the following result we prove that the definition of
standard tangent vector $t \in {\rm T}_{\rm st}M$ is equivalent to the classical one.
\begin{theorem}
    In the category $\CInfty$ the object ${\rm T}_{\rm st}M$ is isomorphic to the usual tangent bundle of $M$
\end{theorem}
\emph{Sketch of the proof:}
We have to prove that ${\rm T}_{\rm st}^m :=\{t \in {\rm T}_{\rm st}M \,|\, t(0)=m \} \simeq {\rm T}_m$ where here
${\rm T}_m := \{ f \in \C^\infty(\R,M) \,|\, f(0)=m  \}/\sim$ is the usual tangent space of $M$ at $m \in M$. Note that
${\rm T}_m \in \CInfty$ because of completeness and co-completeness.\\
Let $d$ be the dimension of $M$. Firstly we prove that
\begin{align*}
    \alpha : \quad [f]_\sim \in {\rm T}_m \quad
    &\mapsto
              \quad \frac{\diff{(\varphi \circ f)}}{\diff t}(0) \in \R^d \\
    \alpha^{-1} : \quad v \in \R^d \quad
    &\mapsto
              \quad [r \mapsto \varphi^{-1}(\varphi m + r \cdot v)]_\sim \in {\rm T}_m
\end{align*}
are arrows of $\CInfty$, where $\varphi : U \freccia \R^d$ is a chart on $m$ with $\varphi(U)=\R^d$.\\
Secondly we prove that
\begin{align*}
    \beta : \quad t \in {\rm T}_{\rm st}^m \quad
    &\mapsto
            \quad \varphi'(t) \in \R^d \\
    \beta^{-1} : \quad v \in \R^d \quad
    &\mapsto
            \quad \ext{[r \mapsto \varphi^{-1}(\varphi m + r \cdot v)]}|_D \in {\rm T}_{\rm st}^m
\end{align*}
are arrows of $\CInfty$. We give some details for $\beta$. If $d \In{U}{\rm T}_{\rm st}^m$ then $d^{\vee} : \bar{U}
\times D \freccia \ext{M}$ in $\ECInfty$. But $\bar{U} \times D = \bar{U} \times \bar{D}=\overline{U \times D}$ hence
$d^{\vee} \In{U \times D} \ext{M}$. Thus we can locally write
$d^{\vee}|_{\mathcal{V}}=\ext{\gamma}(p,-,-)|_{\mathcal{V}}$ where $\mathcal{V}$ is an open neighborhood of $(u,0)$
defined by $A \times B$, $u \in U$ and $\gamma \in \C^\infty(\bar{U} \times A \times B,M)$. But $\mathcal{V} = \ext{(A
\times B)} \cap (U \times D) = (A \cap U) \times D$ because $U \subseteq \R^u$. As in the proof of Theorem \ref{thm:
GeneralizedDerivationFormula} we can prove that
\[
   \beta[d(x)]=\varphi'[d(x)]=\frac{\diff{}}{\diff{t}}\{\varphi[\gamma(p_0,x,t)]\}|_{t=0} \quad \forall x \in A \cap U.
\]
Hence $(d \cdot \beta)|_{A \cap U} \in \C^\infty(A \cap U,\R^d)$ is an ordinary smooth function. Note the importance
to have as $U$ a standard open set in the last passage: this is a strong motivation for the definition we gave of
${\rm T}_{\rm st}M$. \qed
For any object $X \in \ECn$ the multiplication of a tangent vector $t$ by a scalar $r \in \ER$ can be defined simply
``increasing its speed'' by a factor $r$:
\[
    (r \cdot t)(h) := t(r \cdot h).
\]
As we already noted, in the category $\ECn$ we have spaces with singular points too, like algebraic curves with double
points. Because of this reason we cannot always define the sum of tangent vectors, but we need to introduce a class of
objects in which this operation is possible. The following definition simply affirms that in these spaces there always
exists the infinitesimal parallelogram generated by a finite number of given vectors.
\begin{definition}
    Let $X \in \ECn$, then we say that $X$ is \emph{infinitesimally linear} iff for any $k \in \N$ greater
than $1$ and for any $t_i \in {\rm T}_x X$, $i=1,\ldots,k$, there exists one and only one $p : D^k \freccia X$ such
that
\[
  \forall i=1,\ldots,k \pti p(0,\ptind^{i-1},0,h,0,\ldots,0) = t_i(h) \quad \forall h \in D.
\]
\end{definition}
The following theorem gives meaningful examples of infinitesimally linear objects.
\begin{theorem}
    The extension of any manifold $\ext{M}$ is infinitesimally linear. If $M_i \in \Man$ for $i=1,\ldots,s$ then
    \[
        \ext{M_1}^{{\ext{M_2}^{\,\dots}}^{\ext{M_s}}}\simeq \ext{M_1}^{\ext{(M_2 \times \dots \times M_{\scriptstyle s})}}
    \]
    is infinitesimally linear too.
\end{theorem}
\Proof
Given any chart $(U,\varphi)$ on $\st{m}$ we can define the infinitesimal par\-al\-lelo\-gram $p$ as
\begin{equation}
\label{def: InfParallelogram}
    p(h_1,\ldots, h_k) =
        \ext{\varphi}^{-1} \left( \ext{\varphi}(m) + \sum_{i=1}^{k} h_i \cdot \varphi'(t_i) \right).
\end{equation}
If fact if $\tau(h):=p(0,\ptind^{i-1},0,h,0,\ldots,0)$ then $\varphi(\tau(h))=\varphi(m) + h \cdot \varphi'(t_i)$;
this implies that $t(0)=\tau(0)$ and $\varphi'(\tau) = \varphi'(t_i)$, hence $t_i=\tau$. To prove uniqueness consider
that if $p : D^k \freccia \ext{M}$ then $p \In{D^k} \ext{M}$ and we can write $p(h)=\gamma(q,h)$, where $\gamma \in
\C^n(U \times V,M)$ and $q$ is the usual extended parameter. Hence
\[
   \varphi[\gamma(q,0,\ptind^{i-1},0,h,0,\ldots,0)]=\varphi[t_i(h)]=\varphi(m) + h \cdot \varphi'(t_i)
\]
and so
\[
   \varphi[\gamma(q,h)]=\varphi(m) + \sum_{i=1}^k h_i \cdot \varphi'(t_i)
\]
from the first order infinitesimal Taylor formula.\\
Because
\[
    \ext{M_1}^{{\ext{M_2}^{\dots}}^{\ext{M_s}}}
    \simeq
    \ext{M_1}^{\ext{M_2} \times \dots \times \ext{M_{\scriptstyle s}}}
    \simeq
    \ext{M_1}^{\ext{(M_2 \times \dots \times M_{\scriptstyle s})}}
\]
it suffices to prove the conclusion for $s=2$. First of all we note that, because of the previously proved uniqueness,
the definition \ref{def: InfParallelogram} of the infinitesimal parallelogram doesn't depend on the chart $\varphi$ on
$\st{m}$. Now let $t_1, \ldots, t_k$ be $k$ tangent vectors at $f \in \ext{N}^{\ext{M}}$. We shall define their
parallelogram $p : \ext{M} \freccia \ext{N}^{D^k}$patching together smooth functions defined on open subsets, and
using the sheaf property of $\ext{N}^{D^k}$. Indeed for every $m \in \ext{M}$ we can find a chart $(U_m, \varphi_m)$
of $N$ on $\st{f(m)}$ with $\varphi_m(U_m) = \R^n$. Now $m \in V_m := f^{-1}(\ext{U_m}) \in \Top{\ext{M}}$ and for
every $x \in V_m$ we have $t_i^\vee(0,x) = f(x) \in \ext{U_m}$. Hence $t_i^\vee(h,x)\in \ext{U_m}$ for any $h \in D$
by theorem \ref{thm: DandOpenGeneralized}. Therefore we can define
\[
   p^\vee_m(x,h) := \varphi_m^{-1}\left\{ \sum_{i=1}^k \varphi_m[t_i^\vee(h^i,x)] - (k-1) \cdot \varphi_m(fx) \right\}
   \quad \forall x \in V_m, \forall h \in D^k
\]
and we have that $p^\vee_m : (V_m \prec \ext{M}) \times D^k \freccia \ext{N}$ is smooth, because it is composition of
smooth functions. If $x \in V_m \cap V_{m^\prime}$ then $p^\vee_m(x,-) = p^\vee_{m^\prime}(x,-)$, in fact from the
generalized derivation formula $\varphi_m[t_i^\vee(h^i,x)]=\varphi_m(fx)+h^i \cdot \varphi_m^{\prime}[t_i^\vee(-,x)]$
and hence we can write
\begin{equation}
\label{eq: ParallelogramFunctionsSpace}
    p^\vee_m(x,h) = \varphi_m^{-1}\left\{ \varphi_m(fx) + \sum_{i=1}^k h^i \cdot \varphi_m^{\prime}[t_i^\vee(-,x)]
   \right\} \quad \forall x \in V_m, \forall h \in D^k.
\end{equation}
But $(U_m, \varphi_m)$ is a chart on $\st{f(x)}$, so $p^\vee_m(x,-)$ is the infinitesimal parallelogram generated by
the tangent vectors $t_i^\vee(-,x)$ at $f(x)$, and we know that (\ref{eq: ParallelogramFunctionsSpace}) doesn't depend
on $\varphi_m$, so $p_m = p_{m^\prime}$. From (\ref{eq: ParallelogramFunctionsSpace}) is also easy to prove that $p :
D^k \freccia \ext{N}^{\ext{M}}$ verifies the desired properties. Uniqueness follows noting that $p^\vee(m,-)$ is the
infinitesimal parallelogram generated by $t_i^\vee(-,m)$. \qed
If $X$ is infinitesimally linear then we can define the sum of tangent vectors $t_1, t_2 \in {\rm T}_x X$ simply
taking the diagonal of the parallelogram $p$ generated by these vectors
\[
   (t_1 + t_2)(h) := p(h,h) \quad \forall h \in D.
\]
With these operations ${\rm T}_x X$ becomes a $\ER$ module. To prove e.g. that the sum is associative see \cite{Koc,
Lav} for a similar theorem.\\
Vector fields on a generic object $X \in \ECn$ are naturally defined as
\[
    V : X \freccia \hbox{T}X \e{such that} V(m)(0)=m.
\]
In the case of manifolds, $X = \ext{M}$, this implies that $V(m)(0) \in M$ for every $m \in M$, hence from (\ref{thm:
StandardTangentVector})
\[
    V|_{M} : (M \prec \ext{M}) \freccia (\{\ext{f}|_D :\, f \in \Cn(\R,M) \} \prec {\rm T}M).
\]
From this, using the definition of arrow in $\Cn$ and the embedding Theorem \ref{thm: EmbeddingOfSepInExtendedCn}, it
follows that
\[
    V|_{M} : M \freccia {\rm T}_{\rm st}(M) \text{ in } \Cn,
\]
that is the standard notion of vector field on $M$. Vice versa if we have
\[
    W : M \freccia {\rm T}_{\rm st}(M) \text{ in } \Cn
\]
then we can extend it to $\ext{M}$. In fact fix $m \in \ext{M}, h \in D$ and choose a chart $(U,x)$ on $\st{m}$. Then we
can write
\[
    W|_{U} = \sum_{i=1}^{d} A_i \cdot \frac{\partial}{\partial x_i},
\]
with $A_i \in \C^n(U,\R)$. But $m \in \ext{U}$ because $\st{m} \in U$ and hence we can define
\[
    \tilde{W}(m,h) := \sum_{i=1}^{d} \ext{A_i}(m) \cdot \frac{\partial}{\partial x_i}(m)(h).
\]
This definition doesn't depend on the chart $(U,x)$ and, for the sheaf property of $\ext{M}$ provides a $\ECn$
function
\[
    \tilde{W} : \ext{M} \times D \freccia \ext{M} \e{such that} \tilde{W}(m,0) = m
\]
and with $(\tilde{W}^{\wedge})|_{M} = V$.\\
Finally we can easily see that any vector field can equivalently be seen as an infinitesimal transformation of the
space into itself. In fact using cartesian closedness we have
\[
    V \in (X^D)^X \simeq X^{X \times D} \simeq X^{D \times X} \simeq (X^X)^D.
\]
If $W$ corresponds to $V$ in this isomorphism then $W : D \freccia X^X$ and $V(x)(0)=x$ is equivalent to say that
$W(0)=1_X$, that is $W$ is the tangent vector at $1_X$ to the space of transformations $X^X$, that is an infinitesimal
path traced from $1_X$.
\section{A first comparison with other theories of infinitesimals}
It is not easy to clarify in a few rows the relationships between our Infinitesimal Differential Geometry (IDG) and
other, more developed and well established theories of actual infinitesimals. Nevertheless here we want to sketch a
first comparison, and to state some open problems, mostly underlining the conceptual differences instead of the
technical ones, hoping in this way to clarify the foundational and philosophical choices we made.
\subsection{Nonstandard Analysis (NSA)}
As a consequence of the will to have a field which extends the reals, in NSA every non zero infinitesimal is
invertible and so we cannot have nilpotent elements. On the contrary in IDG we aim to obtain a ring as an extension,
and, as a result of our choices, we cannot have non-nilpotent infinitesimals, in particular they cannot be invertible.
In IDG our first aim was to obtain a meaningful theory from the intuitive point of view, to the disadvantage of some
formal property, only partially inherited from the real field. Vice versa any constructions of the hyperreals ${}^*\R$
has, as one of its primary aims, to obtain the inheritance of all the properties of the reals through the transfer
principle. This way of thinking implies that in NSA we want to be free to extend every function $f : \R \freccia \R$
from $\R$ to ${}^*\R$, and that any sequence of standard reals $(x_n)_{n \in \N} \in  \R^{\N}$, even  the more
strange, represents one and only one hyperreal. Of course in this article we followed a completely different way:  to
define $\ER$ we restrict ourselves to nilpotent functions $(x_t)_{t \in \R} \in \Nil \subset \R^{\R}$, and hence we
can only extend locally lipschitzian functions from $\R$ to $\ER$. Obviously we have in mind that in Differential
Geometry we shall work with $\C^n$ functions only. In exchange not every property is transferred to $\ER$, e.g. we
have partial order relations only.\\
In NSA this attention to formally inherit every property of the reals implies that on the one hand we have the
greatest logical strength, but on the other hand we need a higher formal control (some background of Logic is
necessary e.g. to apply the transfer principle) and sometimes we lose the intuitive point of view. E.g. what is the
intuitive meaning and usefulness of $\st{[\sin(I)]} \in \R$, the standard part  of the sine of an infinite number $I
\in {}^*\R$? These, together with very strong but scientifically unjustified cultural reasons, may be some motivations
for the not so high success of NSA in Mathematics, and consequently in its didactics. Anyway NSA is essentially the
only theory of actual infinitesimals with a discrete diffusion and a sufficiently great community of working
mathematicians and published results, even if few of them concern Differential Geometry.\\
Two open problems concerning the relationships between IDG and NSA are the following.
\problem{It is possible to define $\ER$ so as to include the hyperreals of NSA. It suffices to consider sequences
of elements of $\Nil$ and to define
\[
    x \sim y \DIff \left\{n \in \N \,\,|\,\, \Limup{x_n(t) - y_n(t)}=0 \right\} \in \mathcal{U}_\infty.
\]
Where $\mathcal{U}_\infty$ is an ultrafilter which contains the filter of cofinite sets. In this way almost all the
results that we presented here, but not every, can be rightly reformulated. Is it possible to obtain a construction
which follows the ideas presented in this article, but with a good theory of invertible infinitesimals?}
\problem{Our partial order relations are not an order, but we can fix an ultrafilter $\mathcal{U}_0$ which contains the
filter of neighborhoods at $t=0$ and define $\ge$ substituting $\forall^0 t$ in the definition of $\succeq$ with
\[
    \{ t \,|\, x_t \ge y_t + z_t \} \in \mathcal{U}_0,
\]
then we can simply prove that we obtain an order. Modifying in a similar way the equality in $\ER$ is it possible to
prove a general transfer theorem?}
\subsection{Synthetic Differential Geometry (SDG)}
There are many analogies between SDG and IDG, so that sometimes proofs remain almost unchanged. But the differences
are so important that, in spite of the similarities, these theories can be said to describe ``different kind of
infinitesimals''.\\
One of the most important differences is that in IDG we have $h \cdot k = 0$ if $h^2 = k^2 = 0$. This is not the case
in SDG, where infinitesimals $h, k \in \Delta := \{ d \,|\, d^2 = 0 \}$ with $h \cdot k$ not necessarily equal zero,
sometimes play an important role. Note that, as shown in the proof of Schwarz theorem using infinitesimals, to bypass
this difference, sometimes requires completely new ideas. Because of these diversities, in our derivation formula we
are forced to state $\exists ! \,m \in \R$ and not $\exists ! \,m \in \ER$. This is essentially the only important
difference between this formula and the Kock-Lawvere axiom. Indeed to differentiate a generic smooth map $f : \ER
\freccia \ER$ we need ``smooth incremental ratios'' (the analogous of the Fermat-Reyes axiom in SDG; these results of
IDG will be presented in a next work. A first approach to this problem, previous to the introduction of the useful
sheaf property in the definition of $\bar \F$ spaces, can be seen in \cite{Gio2}).

Another point of view of the relationships between these two theories can be introduced starting from a sentence of
\cite{Mo-Re}, pag. 385: \emph{``These structures} [convenient vector spaces] \emph{are in a way simpler than the
sheaves considered in this book, but one should notice that the theory of convenient vector spaces does not include an
attempt to develop an appropriate framework for infinitesimal structures, which is one of the main motivations of our
approach...''}. We want to think that this thought could also be applied to diffeological spaces, and so IDG may be a
possible solution. Indeed models of SDG are not so easy to construct Topos, so that we are almost compelled to work
with the internal language of the Topos itself, that is in intuitionistic logic. If on the one hand this implies that
``all our spaces and functions are smooth'', and so we don't have to prove this after every definition, on the other
hand it requires a more strong formal control of the Mathematics you are doing.\\
Everyone can be in agreement or not with the above cited sentence of \cite{Mo-Re}, or if it is difficult or easy to
learn to work in intuitionistic logic and after to translate the results using Topos models. Anyway we think
undeniable that the formal beauty achieved by SDG can with difficulty be reached using a theory in classical logic. It
suffices to say, as a simple example, that to prove the infinitesimal linearity of $M^N$ (starting from $M$, $N$
generic infinitesimally linear spaces), it suffices to fix $n \in N$, to note that $t_i(-,n)$ are tangent vectors at
$f(n)$, to consider their parallelogram $p(-,n)$, and automatically, thanks to the use of intuitionistic logic, $p$ is
smooth without any need to use directly the sheaf property to prove it.\\
On the other hand if we need a partition of unity, we are forced to assume a suitable axiom for the existence of bump
functions (whose definition, in the models, necessarily uses the law of the excluded middle).

From the intuitive, classical, point of view, it is a little strange that we don't have ``examples'' of infinitesimals
in SDG (it is only possible to prove that $\neg\neg \exists d \in \Delta$), so that, e.g., we cannot construct a
physical theory containing a fixed infinitesimal parameter; moreover any $d \in \Delta$ is at the same time $d \le 0$
and $d \ge 0$; finally the definition of the Lie brackets using $h \cdot k$ for $h, k \in \Delta$ is very far to the
usual definitions given on manifolds.
\problem{Is it possible to construct a theory of nilpotent infinitesimals useful for several construction in Differential
Geometry and with:
\begin{itemize}
    \item meaningful and useful examples of first order infinitesimals $h^2=k^2=0$ with $h \cdot k \ne 0$;
    \item models simpler than Topos models of SDG so that classical logic suffices to work in it;
    \item $\exists ! \, m \in \ER$ in the derivation formula?
\end{itemize} }
\subsection{Weil functors (WF)}
Weil functors (see \cite{Ko-Mi-Sl}) represent, as far as we know, the only way to introduce some kind of useful
infinitesimal method without the need to possess a non-trivial background in mathematical logic. They don't arrive to
the construction of a whole ``infinitesimal universe'' like in IDG or in the previously cited theories, but to define
functors $T_A  : \ManInfty \freccia \ManInfty$, related to the geometrical constructions we are interested in and
starting from a \emph{Weil algebra} $A = \R \cdot 1 \oplus N$ ($N$ is a finite dimensional ideal of nilpotent
elements). The flexibility of its input $A$ gives a corresponding flexibility to the construction of these functors.
But, generally speaking, if one change the geometrical problem, one has also to change the algebra $A$ and so the
corresponding functor $T_A$. E.g. if $A = \R[x]/\langle x^2 \rangle$, then $T_A$ is the ordinary tangent bundle
functor, whereas if $B = R[x,y]/\langle x^2, y^2 \rangle$, then $T_B = T_A \circ T_A$ is the second tangent bundle.
Note that $x$, $y \in B$ verify $x^2=y^2=0$ but $x \cdot y \ne 0$. This provides us with the first difference between
WF and IDG. In fact although $\ER = \R \cdot 1 \oplus I_0$ and $\dim_{\R}I_0 = \infty$, using the infinitesimals of
$\ER$ we can generate a large family of Weil algebras (e.g. any $A = \R \cdot 1 \oplus N \subset \R \cdot 1 \oplus
D_k$ which represents $k$th order infinitesimal Taylor formula) but not every algebra can be generated in this way,
e.g. the previous $B$. But using exponential objects of $\CInfty$ and $\ECInfty$ we can give a simple infinitesimal
representation of a large class of WF. We will use the common multi-index notations, e.g. if $\alpha = (\alpha^1,
\ldots, \alpha^n) \in \N^n$, $h = (h_1, \ldots, h_n)\in \ER^n$, then $h^\alpha = h_1^{\alpha^1}\cdot \ldots \cdot
h_n^{\alpha^n} \in \ER$. For $\alpha_1, \ldots , \alpha_c \in \N^n$, $c \ge n$, let
\[
   D_\alpha := \left\{ h \in \ER^n \,|\, h^{\alpha_i} \in D \ \ \forall i=1, \ldots ,c \right\}.
\]
Using this notation we will always suppose that $\alpha$ verifies
\begin{align*}
    {}& \alpha_1 = (k_1,0,\ldots,0) \\
    {}& \alpha_2 = (0,k_2,0,\ldots,0) \\
    {}& \ldots\\
    {}& \alpha_n = (0,\ldots,0,k_n)\\
    {}& \alpha_i^j \le k_j \quad \forall i=n+1, \ldots, c.
\end{align*}
Hence $D_\alpha \subseteq D_{k_1} \times \ldots \times D_{k_n}$. E.g. if $\alpha_1 =(3,0)$, $\alpha_2 =(0,2)$ and
$\alpha_3 =(1,1)$, then $D_\alpha = \{(h,k) \in D_3 \times D_2 \,|\, h \cdot k \in D\}$. To any infinitesimal object
$D_\alpha$ is associated a corresponding Taylor formula: let $f = \ext{g}|_{D_\alpha}$, with $g \in
\C^\infty(\R^n,\R)$, then
\[
   f(h) = \sum_{
                \substack{r \in \N^n\\
                          \exists i: \, r \le \alpha_i
                         }
               }
                \frac{h^r}{r!} \cdot m_r \quad\forall h \in D_\alpha.
\]
Coefficients $m_r = \frac{\partial^r g}{\partial x^r}(0)\in \R$ are uniquely determined by this formula. Here $r \le
\alpha_i$ means $r^j \le \alpha_i^j$ for every $j=1, \ldots, n$. We can therefore proceed generalizing the definition
\ref{def: StndTangentFunctor} of standard tangent functor.
\begin{definition}
We call $M^{D_\alpha}$ the $\CInfty$ object with support set
\[
    \left| M^{D_\alpha} \right| := \{\ext{f}|_{D_\alpha} \,:\, f \in \CInfty(\R^n,M) \},
\]
and with generalized elements of type $U$ (open in $\R^u$)
\[
    d \In{U} M^{D_\alpha} \DIff d: U \freccia \left| M^{D_\alpha} \right|
          \e{and} d \cdot i \In{\bar{U}} \ext{M}^{D_\alpha},
\]
where $i : \left| M^{D_\alpha} \right| \hookrightarrow \ext{M}^{D_\alpha}$ is the inclusion.
\end{definition}
We can extend this definition to the arrows of $\ManInfty$ with $f^{D_\alpha}(t) := t \cdot f \in N^{D_\alpha}$, where
$t \in M^{D_\alpha}$ and $f \in \ManInfty(M,N)$. With these definitions we obtain a product preserving functor
$(-)^{D_\alpha}: \ManInfty \freccia \ManInfty$. Finally we have a natural transformation $e_0 : (-)^{D_\alpha}
\freccia 1_{\ManInfty}$ defined by evaluation at $0 \in \R^n$:  $e_0(M)(t) := t(0)$. The functor $(-)^{D_\alpha}$ and
the natural transformation $e_0$ verify the ``locality condition'' of theorem 1.36.1 in \cite{Ko-Mi-Sl}: if $U$ is
open in $M$ and $i : U \hookrightarrow M$ is the inclusion, then $U^{D_\alpha} = e_0(M)^{-1}(U)$ and $i^{D_\alpha}$ is
the inclusion of $U^{D_\alpha}$ in $M^{D_\alpha}$. We can thus apply the above cited theorem to obtain that
$(-)^{D_\alpha}$ is a Weil functor, whose algebra is
\[
   Al\left( (-)^{D_\alpha} \right) = \R^{D_\alpha} \simeq \R[x_1, \ldots, x_n]/ \langle x_1^{\beta^1} \cdot \ldots
   \cdot x_n^{\beta^n} \rangle_{i \in I}.
\]
Where $I := \{ \beta \in \N^n \,|\, \exists i \exists j : \, \beta^j > \alpha_i^j \}$.\\
Not every Weil functor has this simple infinitesimal representation. E.g. the second tangent bundle $(-)^D \circ
(-)^D$ is not of type $(-)^{D_\alpha}$; indeed it is easy to prove that the only possible candidate could be $D_\alpha
= D \times D$, but $(\R^D)^D$ is a four dimensional manifold, whereas $\R^{D \times D}$ has dimension three. We don't
have this kind of problems with the functor $(-)^{D_\alpha} = \ECInfty(D_\alpha,-): \ECInfty \freccia \ECInfty$ which
generalizes the previous one as well as ${\rm T}M = \ext{M}^D$ generalizes the standard tangent functor. In fact
because of cartesian closedness we have
\[
   \left( X^{D_\alpha} \right)^{D_\beta} \simeq X^{D_\alpha \times D_\beta}
\]
and $D_\alpha \times D_\beta$ is again of type $D_\alpha$.\\
Weil functors has another more general, but less simple, infinitesimal representation using exponential objects. We
sketch here the case of the second tangent bundle for $M=\R$ only. Let
\begin{align*}
    {}&D(10) := \left\{ (h,k) \in D_2^2 \,|\, h^2=0=k \right\} \\
    {}&D(01) := \left\{ (h,k) \in D_2^2 \,|\, h=0=k^2 \right\} \\
    {}&D(11) := \left\{ (h,k) \in D_2^2 \,|\, h^2 \cdot k=0=h \cdot k^2 \right\} \\
    {}&\bar D := D(10) \times D(01) \times D(11).
\end{align*}
Now we consider the incremental differences corresponding to these objects, that is
\begin{align*}
    f_{10}[-] &: (h,k) \in D(10) \mapsto f(h,k) - f(0) \in \ER \\
    f_{01}[-] &: (h,k) \in D(01) \mapsto f(h,k) - f(0) \in \ER \\
    f_{11}[-] &: (h,k) \in D(11) \mapsto f(h,k) - f(h,0) - f(0,k) + f(0) \in \ER.
\end{align*}
Finally let
\[
   \R^{\bar D} := \left\{ (f(0), p_1 \cdot f_{10}[-], p_2 \cdot f_{01}[-], p3 \cdot f_{11}[-])
                         \,|\, f \in \C^{\infty}(\R^2,\R)
                 \right\}
\]
where $p_i$ are projections, e.g. $p_1 : \bar D \freccia D(10)$. If $g \in \R^{\bar D}$ then
\[
   g(h) = (f(0), h_1 \cdot \partial_1 f(0), h_4 \cdot \partial_2 f(0), h_5 \cdot h_6 \cdot \partial_{12} f(0))
   \quad \forall h \in \bar D
\]
and $\R^{\bar D}$ is an algebra isomorphic to $(\R^D)^D$. To generalize this representation to a generic manifold $M
\in \ManInfty$ we have to use infinitesimal incremental differences of functions $f \in \C^\infty(\R^n, M)$. This
could also be performed in an intrinsic way using affine structures definable on the infinitesimals neighborhood of
any point, which will be presented in a next work.
\section*{Acknowledgment}
The author wish to thank Prof. Sergio Albeverio for his great encouragement and support.
\bibliographystyle{plain}
\nocite{*}
\bibliography{giordano}
\end{document}